\documentclass[11pt]{amsart}

\usepackage{amssymb,mathtools}
\usepackage{bbm}
\usepackage[normalem]{ulem}

\newcommand{\beq}{\begin{equation}}
\newcommand{\eeq}{\end{equation}}
\numberwithin{equation}{section}

\usepackage{enumitem,color}

\newtheorem{thm}{Theorem}[section]

\newtheorem {lemma}[thm]{Lemma}
\newtheorem {cor}[thm] {Corollary}

\newtheorem{thmA}{Theorem}

\theoremstyle{definition}

\newtheorem{eg}{Example}[section] 
\newtheorem{rmk}{Remark}[section]

\newcommand{\du}{\,\mathrm{d}\mkern0.6mu u}
\newcommand{\dv}{\,\mathrm{d}\mkern0.6mu v}

\newcommand{\Fock}[2]{\mathcal{F}^{#1}_{#2}}

\newcommand{\F}{\Fock{p}{\alpha}}
\newcommand{\Fh}{\Fock{2}{\alpha}}
\newcommand{\Fq}{\Fock{q}{\alpha}}
\newcommand{\Finfty}{\Fock{\infty}{\alpha}}
\newcommand{\Fzero}{\Fock{\infty}{\alpha,0}}
\newcommand{\Fone}{\Fock{1}{\alpha}}

\newcommand{\C}{\mathbb{C}}
\newcommand{\R}{\mathbb{R}}
\newcommand{\N}{\mathbb{N}}
\newcommand{\dm}{\,\mathrm{d}\mkern0.6mu m}
\newcommand{\comp}{C_\varphi}
\newcommand{\wcomp}{W_{\psi,\varphi}}
\newcommand{\mcomp}{M(\psi,\varphi)}
\newcommand{\mzcomp}{M_z(\psi,\varphi)}
\newcommand{\twcomp}{\widetilde W_{\psi,\varphi}}

\newcommand{\wbar}{\overline{w}}

\newcommand{\abs}[1]{\left| #1 \right|}
\newcommand{\norm}[1]{{\left\|#1\right\|}}

\newcommand{\ip}[2]{\langle #1,\, #2 \rangle}  

\newcommand{\spn}[1]{\mathrm{span}\!\left\lbrace #1 \right\rbrace} 

\usepackage[numbers]{natbib}

\title[Weighted composition operators on Fock spaces]{Weighted composition operators on Fock spaces
	and their dynamics} 
\author[T. Carroll]{Tom Carroll}                                
\address{School of Mathematical Sciences, University College Cork, Ireland}                                                                    
\email{t.carroll@ucc.ie}                                                                          
\author[C. Gilmore]{Clifford Gilmore}                                
\address{School of Mathematical Sciences, University College Cork, Ireland}                                                              
\email{clifford.gilmore@ucc.ie}                                                                                                   
\date{}                                       
                                                   

\thanks{C.~Gilmore was supported by the Irish Research Council via a Government of Ireland Postdoctoral Fellowship.}

\begin{document}                                                

\begin{abstract}
Bounded weighted composition operators, as well as compact weighted composition  
operators, on Fock spaces have been characterised. 
This characterisation is refined to the extent that the 
question of whether weighted composition operators on the Fock space 
can be supercyclic is answered in the negative. 
\end{abstract}                

\maketitle  

\section{Introduction}\label{sec:intro}

The study of composition operators $C_\varphi f = f \circ \varphi$ acting on spaces of analytic functions has been an active area of research since the pioneering work of Nordgren~\cite{Nor68} and Kamowitz~\cite{Kam75, Kam78}.  Although this class of operators is usually tractable, there exists a rich literature that reveals many interesting properties in a variety of settings, including in the Hardy, Bergman and Dirichlet spaces~\cite{NRW87}, \cite{Cow83}, \cite{Sha87}, \cite{GS16}.  

A natural generalisation of the family of composition operators are the weighted composition operators $\wcomp f = \psi \cdot (f\circ \varphi)$, which have recently attracted much research interest~\cite{CH01}, \cite{AL04}, \cite{Gun07}, \cite{Gun11},  \cite{BLW12}, \cite{HLNS13}, \cite{TK19}.  In particular, the bounded and compact weighted composition operators acting on Fock spaces of entire functions were characterised by Ueki~\cite{Uek07, Uek10}, Le~\cite{Le14}, Hai and Khoi \cite{HK16}, 
and Tien and Khoi \cite{TK19}.  
Analogous results for unweighted composition operators acting on the classical Fock space were previously identified by Carswell et al.~\cite{CMS03}.

Hitherto, characterisations of the boundedness and compactness of $\wcomp$ in the Fock space setting 
have required that the symbol $\varphi$ and the multiplier $\psi$ satisfy certain uniform conditions.  
In Section \ref{sec:lessthan1} of this article we study $\wcomp$ from the perspective of the order and type of the multiplier $\psi$.  
This approach allows us to obtain more explicit boundedness and compactness conditions (Theorem~\ref{thm1}). 
In particular, we give a precise description of the zero-free multipliers $\psi$ 
that admit bounded or compact operators
(Theorem~\ref{thm2}).

The second theme of this article is to study the linear dynamical properties of weighted composition operators  in the Fock space setting.  
Linear dynamics has been an active branch of operator theory since the early 1990s, and composition operators acting on various function spaces have been studied from this perspective in \cite{BS90}, \cite{BM95}, \cite{BS97}, \cite{GGMR04}, \cite{BD11}, \cite{BBMP15}, \cite{Bes13}, \cite{BGJJ16a}, \cite{Tie20a} and \cite{Tie20b}.  
Analogous questions for weighted composition operators were investigated in \cite{YR07}, \cite{KHK10}, \cite{Rez11}, \cite{Bes14}, \cite{BGJJ16b} and \cite{Bel20}. 
Comprehensive introductions to linear dynamics can be found in the monographs by Bayart and Matheron~\cite{BM09}, and Grosse-Erdmann and Peris~\cite{GEP11}. A survey of some of the recent advances in the area can be found in \cite{Gil20}.

To investigate the linear dynamical properties of weighted composition operators acting on Fock spaces, 
we use the results from Section \ref{sec:lessthan1} to  give explicit asymptotics for the iterates of $\wcomp$ in this setting (Theorems~\ref{lambdacomplex}, \ref{propreal} and  \ref{prop:realinfinity}). As an application of these results, we show that $\wcomp$ acting on Fock spaces cannot be supercyclic (Theorem~\ref{thm3}).  
This illustrates that weighted composition operators have contrasting dynamical behaviours depending on the setting.

\section{Background}\label{sec:Fock}

\subsection{Fock spaces}

\noindent In this section we recall the pertinent notation and results, taken from Zhu~\cite{Zhu12}, which are needed in the sequel.

For $\alpha$ positive and $1\leq p <\infty$, the Fock space $\F$ is composed of 
entire functions $f$ for which the norm 
\beq\label{1.1}
\norm{f}_{p,\alpha}
=  \left( \frac{\alpha p}{2\pi}\int_\C \abs{ f(z)}^p\,e^{-\tfrac{\alpha p}{2}\abs{ z}^2}\, \dm(z) \right)^{1/p}
\eeq
is finite. 
Here $\dm$ is area measure in the complex plane $\C$. 

For $\alpha$ positive and $p=\infty$, the Fock space $\Finfty$ is composed of 
entire functions $f$ for which the norm 
\beq\label{Finftynorm}
\norm{f}_{\infty,\alpha} = \sup \left\lbrace \abs{f(z)} 
e^{-\tfrac{\alpha}{2}\abs{z}^2} : z \in \C  \right\rbrace 
\eeq
is finite. 
Henceforth, the dependence of the norm on $\alpha$ and on $p$ is suppressed in the notation.

We let  $L_\alpha^p$ denote the space of Lebesgue measurable functions $f$ on $\C$ such that the function $f(z)e^{-\alpha\abs{z}^2/2} \in L^p(\C, \dm)$.
Thus an entire function $f$  belongs to the Fock space $\F$
when $f \in L_\alpha^p$. It is well known that $\F$ is a Banach space  for $1\leq p\leq \infty$.

We will also consider the closed subspace $\Fzero$ of $\Finfty$, which consists of the entire functions $f$ such that
\begin{equation*}
\lim_{z\to \infty} f(z) e^{- \frac{\alpha}{2} \abs{z}^2} = 0.
\end{equation*}
The space $\Fzero$ is the closure in $\Finfty$ of the set of polynomials.  So in particular we have that $\Fzero$ is separable, while $\Finfty$ is non-separable (cf.~\cite[p.~39]{Zhu12}).

Point evaluations are bounded linear functionals on $\F$ for $1 \leq p \leq \infty$. This follows immediately from the following pointwise estimate: for $f \in \F$ 
\begin{equation} \label{ineq:ptEvalFinite}
\abs{f(z)} \leq \norm{f} e^{\alpha \abs{z}^2/2}
\end{equation}
for all $z \in \C$ and this estimate is sharp (cf.~\cite[Corollary 2.8]{Zhu12}).

For $1 \leq p < \infty$, the dual space of $\F$ can be identified with $\Fq$ via the integral pairing  
\[
\ip{f}{g} = \frac{\alpha}{\pi}\int_\C f(w)\,\overline{g(w)}
\,e^{-\alpha \vert w \vert^2}\, \dm(w)
\]
for $g \in \Fq$ and $f \in \F$ (cf.~\cite[Corollary 2.25]{Zhu12}).
Here $p$ and $q$ are the usual conjugate exponents with $1/p + 1/q = 1$.
In the case $\Fzero$, the dual space can be identified with $\Fone$ (cf.~\cite[Theorem 2.26]{Zhu12}).

For $z\in \C$ and $\alpha$ positive,  
we define the kernel function as
\begin{equation*}
k_z(w)  = e^{\alpha\overline{z}w}, \quad   w\in \C. 
\end{equation*}
It has norm $\Vert k_z \Vert = e^{\alpha\vert z\vert^2/2}$ and belongs to $\F$ for $1\leq p\leq \infty$.

It turns out that point evaluation on $\F$,  for $1 \leq p \leq \infty$, 
is realised through the kernel functions $k_z$. 
For the convenience of the reader we briefly outline why this holds. 
Following Zhu~\cite[Section 2.2]{Zhu12}, for $\alpha>0$ and $1\leq p\leq \infty$, we define the integral operator $P_\alpha$ on $L_\alpha^p$ by 
\begin{equation*}
P_\alpha f(z) =  \frac{\alpha}{\pi} \int_{\C} e^{\alpha z \wbar } f(w) e^{-\alpha \abs{w}^2}\dm(w), 
\quad f \in L^p_\alpha
\end{equation*}	
and we note that $P_\alpha \colon L^p_\alpha \to \F$ is a bounded operator \cite[Corollary~2.22]{Zhu12}.	
Moreover, for $f\in\F$, we have that $P_\alpha f = f$. To see why this holds, we first note that since $\F\subseteq\Finfty$, it is enough to check this for $f \in \Finfty$. 
If $h\in \Finfty$, then $\vert h(\zeta) \vert \leq Ce^{\alpha\vert \zeta\vert^2/2}$, $\zeta\in\C$ and some constant $C$. 
Integrating over circles and using the mean value property we obtain
\[
h(0) = \frac{\alpha}{\pi}\int_\C h(\zeta)e^{-\alpha \vert \zeta \vert^2}\,\dm(\zeta).
\]
For $f\in \Finfty$ and $z \in \C$, we apply the previous formula to 
$h(\zeta) = e^{-\alpha \overline{z}\zeta} f(\zeta+z)$, which also belongs to $\Finfty$,
to obtain that
\begin{align} \label{integralRep}
f(z) 
& = \frac{\alpha}{\pi}\int_\C f(\zeta+z)  e^{-\alpha( \abs{\zeta}^2 + \overline{z}\zeta)} \, \dm(\zeta)
\nonumber\\
& = \frac{\alpha}{\pi}\int_\C f(w)\,\overline{k_z(w)}\,e^{-\alpha \abs{w}^2}\, \dm(w).
\end{align}   
So it follows that $P_\alpha f = f$ for all $f \in \Finfty$. 
That $P_\alpha \colon L_\alpha^p \to \F$, $1 \leq p \leq \infty$, is a bounded projection can also be deduced from a more general theorem in the scholium by Janson et al.~\cite[Theorem 7.1]{JPR87}.

In the case $p=2$, the Fock space $\Fh$ is 
a reproducing kernel Hilbert space with reproducing kernel at $z$ given by $k_z$.
The normalised point evaluations act weakly on functions $f$ in the Fock space 
$\F$ for $p<\infty$, that is
\beq\label{weaklyto0}
\abs{ \ip{f}{k_z/ \norm{ k_z}}} 
= \abs{f(z)}\,e^{-\alpha \abs{ z}^2/2} \to 0, \quad  \textrm{ as }	\abs{z} \to \infty.
\eeq

\begin{rmk} \label{rmk:QuadraticsInFock}
	The entire function $f(z) = e^{az + b z^2}$ belongs to the Fock space $\F$, 
	$1\leq p < \infty$, or to $\Fzero$,
	if and only if $\abs{b} < \alpha/2$. 
	The function $f(z) = e^{b z^2}$ belongs to $\Finfty$ if
	and only if $\vert b \vert \leq \alpha/2$.
	The fact that $f(z) = e^{\alpha z^2/2}$ belongs to $\Finfty$, 
	but not to $\Fzero$ or to $\F$ for $p<\infty$, 
	will be of some significance  in our later observations.
\end{rmk}

We note that some articles only consider the Fock spaces corresponding to the case $\alpha=1$. 
The extension to general $\alpha$ is straightforward
and when we recall results in this generality we do so without additional comment. 

\subsection{Weighted composition operators acting on  Fock spaces}\label{sec:lit}

Carswell et al.~\cite{CMS03} characterised the bounded 
(unweighted) composition operators acting on the Fock space $\Fh$ of functions holomorphic on $\C^n$. 
For our present purposes, we state their result for the base case  $n=1$. 
The composition operator $\comp$ with symbol $\varphi$ is defined as
$\comp f = f\circ\varphi$, for $f \in \F$.  
The interested reader can find comprehensive introductions to the topic of composition operators in the monographs by Shapiro~\cite{Sha93}, and Cowen and MacCluer~\cite{CM95}.  

The pertinent result from \cite{CMS03} is as follows.

\begin{thmA}[Carswell, MacCluer and Schuster~\cite{CMS03}]\label{thmA}The 
	composition operator $\comp$ is bounded on the Fock space $\Fh$ only 
	when $\varphi(z) = a + \lambda z$ with $\vert \lambda \vert \leq 1$. 
	
	Conversely, suppose that  $\varphi(z) = a+\lambda z$.
	If $\vert \lambda \vert =1$  then 
	$\comp$ is  bounded  on $\Fh$ if and only if $a=0$. 
	If $\vert \lambda \vert < 1$ then $\comp$ acts compactly on $\Fh$. 
\end{thmA}

Following the initial work of Ueki~\cite{Uek07}, Le~\cite{Le14} 
generalised the results of \cite{CMS03} to weighted 
composition operators acting on the Fock space $\Fh$. 
The extension to $\F$, for $1\leq p \leq \infty$,
was carried out in part by Hai and Khoi \cite{HK16}, 
while Tien and Khoi \cite{TK19} characterised boundedness and compactness of 
weighted composition operators between different Fock spaces and for the full range of $p$.

The weighted composition operator $\wcomp$, with symbol 
$\varphi$ and multiplier $\psi$, is defined as
\[
\wcomp f = \psi \cdot (f\circ \varphi)
\]  
for $f \in \F$.  We assume throughout that the multiplier $\psi$ is not identically zero. 

A necessary condition for the boundedness of $\wcomp$ acting on $\F$ is deduced 
from the following elementary observation. 
Since the constant function $1 \in \F$, it follows for bounded $\wcomp$ on $\F$ 
that $\psi$ itself must belong to $\F$.

For $z \in \C$, we know from \eqref{ineq:ptEvalFinite}  that for $1\leq p \leq \infty$ the point evaluation 
$k_z$ is a continuous linear functional on $\F$.  
If we consider a bounded weighted composition operator $\wcomp$ acting on $\F$, 
for $1\leq p < \infty$, then by Fock space duality the action of the 
(Banach space) adjoint $\wcomp^*$ on   $k_z \in \Fq$ (where $1/p + 1/q = 1$),  is as follows: for any $f \in \F$ 
\begin{align*} 
\left( \wcomp^* k_z \right)(f) &= \ip{\wcomp^* k_z}{f} \\
&= \ip{k_z}{ \wcomp f} \\
&= ( \wcomp f)(z) \\
&= \psi(z) f(\varphi(z)) \\
&= \psi(z) k_{\varphi(z)}(f).
\end{align*}

A similar calculation holds for bounded operators $\wcomp \colon \Fzero \to \Fzero$, where the adjoint $\wcomp^*$ acts on $k_z \in \Fone$. Thus  the action of the adjoint $\wcomp^*$ 
on the pertinent kernel function $k_z$ is given by
\begin{equation}  \label{wstar}
\wcomp^* k_z = \psi(z) k_{\varphi(z)}.
\end{equation}
As a consequence,
\beq\label{2.1}
\frac{\norm{\wcomp^* k_z}^2}{\norm{k_z}^2}
= \frac{\abs{ \psi(z)}^2 \norm{ k_{\varphi(z)}}^2}{\norm{k_z}^2}
= \abs{\psi(z)}^2\,  e^{ \alpha \left( \abs{\varphi(z)}^2 - \abs{z}^2 \right) }
\eeq
is uniformly bounded in $z$. 

For $z \in \C$ and  entire functions $\psi$ and $\varphi$,  we set 
\begin{align*} 
\mzcomp &\coloneqq \vert \psi(z) \vert^2\, e^{ \alpha \left( \vert \varphi(z) \vert^2 - \vert z\vert^2 \right) }
\shortintertext{and}
\mcomp &\coloneqq \sup\left\{ \mzcomp : z\in\C \right\}.   
\end{align*}

We combine  below in Theorem \ref{thmB} the results from \cite{Le14} and \cite{HK16} 
on the characterisation of the bounded weighted composition operators acting on Fock spaces.

\begin{thmA}[Le~\cite{Le14}, Hai and Khoi~\cite{HK16}] \label{thmB}
	Let $\wcomp$ be the weighted composition operator  acting on the Fock space $\F$, $1 \leq p \leq \infty$, or on $\Fzero$. Then $\wcomp$ is bounded  if and only if both 
	$\psi \in \F$, respectively $\Fzero$,
	and $\mcomp  < \infty$.  
\end{thmA}

Although not explicitly stated in the literature, we remark that Theorem \ref{thmB} also holds in the setting of $\Fzero$. 
To see why this is the case, note that since both $1\in \Fzero$ and 
$k_z \in \Fone = \left( \Fzero \right)^*$, 
the conditions  that $\psi$ belongs to $\Fzero$ and that $\mcomp  < \infty$ are again 
necessary for the boundedness of $\wcomp$ on $\Fzero$. 
Assuming these conditions, 
$\wcomp$ is bounded on $\Finfty$. 
{To see that} $\Fzero$ is invariant under $\wcomp$, take 
$f\in \Finfty$ 
and $z\in \C$, so that 
\begin{align}\label{WonFinfty}
\vert (\wcomp f)(z)\vert \,e^{-\tfrac{\alpha}{2} \vert z \vert^2} & =
\vert \psi(z)\vert\, \vert f(\varphi(z)) \vert\, e^{-\tfrac{\alpha}{2} \vert z \vert^2} \nonumber\\
& = \sqrt{ \mzcomp } \,e^{-\tfrac{\alpha}{2} \vert \varphi(z) \vert^2}\,\vert 
f(\varphi(z)) \vert \nonumber\\
& \leq \sqrt{\mcomp}\, \vert f(w) \vert\, \,e^{-\tfrac{\alpha}{2} \vert w \vert^2},
\end{align}
where we let $w=\varphi(z)$.  
Taking the supremum over $z$ we see that $\wcomp$ is bounded on $\Finfty$. Moreover,
when $\lambda \neq 0$ we have that $w = a + \lambda z \to \infty$ as $z\to \infty$, 
and if $f\in \Fzero$, it follows from \eqref{WonFinfty} that $\wcomp f \in \Fzero$, 
so that $\Fzero$ is invariant under $\wcomp$.

We note the condition that $\mcomp$ be finite does not depend on $p$.
In \cite[Proposition 2.1]{Le14} it was shown that, if 
$\mcomp<\infty$ for entire functions $\psi$ and $\varphi$ with $\psi$ not identically zero, 
then $\varphi$ is of the form
\begin{equation} \label{defn:formOfPhi}
\varphi(z) = a + \lambda z 
\end{equation}
where $\abs{\lambda} \leq 1$ and $a \in \C$.
In the case $\vert \lambda \vert = 1$, it was also shown in \cite[Proposition 2.1]{Le14} 
that the multiplier $\psi$ is of the form
\beq\label{defn:formOfPsi}
\psi(z) = \psi(0) e^{- \alpha\overline{a}\lambda z} = \psi(0) k_{-a\overline{\lambda}}(z),
\eeq
in which case $\mzcomp$ is constant and positive.

\begin{rmk} \label{rmk:lambdaEqZero}
	As mentioned in \cite[Corollary~3.2]{TK19}, 
	the case $\lambda=0$, which gives that $\varphi$ is constant, is both trivial and anomalous. 
	The weighted composition operator $\wcomp$ reduces to $(\wcomp f)(z) = \psi(z) f(a)$. 
	The condition for boundedness is simply that $\psi$ belongs to the space in question and in this case the operator is of finite rank and hence compact. 
	We will thus generally assume that $0 < \vert \lambda \vert \leq 1$. 
\end{rmk}

The compact weighted composition operators acting on $\Finfty$ were characterised by Ueki~\cite{Uek10}, and their characterisation in $\Fh$ was given in \cite{Le14}.  These characterisations were extended in \cite{HK16} to the Fock spaces $\F$ for $1 < p < \infty$, and the case $p=1$ can be shown using the approach from \cite{TK19}. 
The various results on compactness are combined into the following theorem.

\begin{thmA}[Le~\cite{Le14}, Ukei~\cite{Uek10}, Hai and Khoi~\cite{HK16}, Tien and Khoi~\cite{TK19}]\label{thmC}
	Let $\varphi(z) = a + \lambda z$, with $0 < \vert \lambda\vert \leq 1$, and let $\wcomp$ be the 
	bounded weighted composition operator acting on $\F$, $1\leq p \leq \infty$, or on $\Fzero$. 
	Then $\wcomp$ is compact if and only if 
	\beq\label{2.2}
	\mzcomp \longrightarrow 0, \textrm{ as } \vert z\vert \to \infty.
	\eeq
\end{thmA}

Although not explicitly stated in the above references, Theorem \ref{thmC} holds 
for $\wcomp$ acting on the space $\Fzero$. To see why this is the case, 
suppose that $\psi \in \Fzero$ and that \eqref{2.2} 
holds. Then $\wcomp$ is compact on $\Finfty$ by Theorem~\ref{thmC}, and since $\mcomp< \infty$ we have that $\Fzero$ is invariant under $\wcomp$. 
Hence, $\wcomp$ acts compactly on $\Fzero$. 

Conversely, suppose that $\wcomp$ acts compactly on $\Fzero$. 
Let $\{z_n\}_n \subset \C$ be a sequence tending to infinity such that the sequence
$\{w_n\}_n = \{\varphi(z_n)\}_n$ also tends to infinity. 
The sequence $\{k_{w_n} / \Vert k_{w_n} \Vert\}_n$ converges weakly to 0 in 
$\Fzero$. In fact, for $f \in \Fone$, 
\[
\abs{ \langle f , k_{w_n} / \norm{k_{w_n}} \rangle}
= \abs{ f(w_n) } \, e^{-\tfrac{\alpha}{2} \abs{w_n}^2} 
\longrightarrow  0 \mbox{ as } n \to \infty.
\]
Hence the sequence $\left\lbrace \wcomp k_{w_n} / \Vert k_{w_n} \Vert \right\rbrace_n$ converges to 0 in $\Fzero$,
so that 
\[
\sup_{z \in \C} \left\{ \left( \vert (\wcomp k_{w_n})(z) \vert\, 
e^{-\tfrac{\alpha}{2}\vert w_n \vert^2}\right)\,  
e^{-\tfrac{\alpha}{2}\vert z \vert^2}   \right\}
\longrightarrow 0, \textrm{ as } n \to \infty.
\]
Choose $z=z_n$ so that 
\[
(\wcomp k_{w_n})(z_n) = \psi(z_n)\, k_{w_n}(\varphi(z_n)) 
= \psi(z_n) k_{w_n}(w_n) = \psi(z_n) e^{\alpha \vert w_n\vert^2},
\]
from which we conclude that  
\begin{align*}
\left( \vert\psi(z_n)\vert\, e^{\alpha \vert w_n\vert^2} \right)
e^{-\tfrac{\alpha}{2} \vert w_n\vert^2}e^{- \tfrac{\alpha}{2} \vert z_n\vert^2}
& = \vert \psi(z_n)\vert\, e^{\tfrac{\alpha}{2}(\vert \varphi(z_n) \vert^2 
	- \vert z_n\vert^2)} \\
& = M_{z_n}(\psi,\varphi)  \longrightarrow  0  \mbox{ as } n \to \infty.
\end{align*}
Hence $\mzcomp \to 0$ as $z \to \infty$ and \eqref{2.2} holds.

We know from \eqref{defn:formOfPsi} that $\mzcomp$ is constant and positive when
$\varphi(z) = a + \lambda z$ with $\vert \lambda \vert=1$.  Hence $\wcomp$ cannot be compact for $\vert \lambda\vert=1$. 
The converse, that $\vert \lambda \vert < 1$ implies compactness of $\wcomp$, 
does not hold in general. Indeed, we will see that this case requires a more careful analysis.

\section{Boundedness and compactness of $\wcomp$ when $\vert\lambda\vert < 1$}\label{sec:lessthan1}

We consider here symbols $\varphi(z)=a+\lambda z$, for $\vert\lambda\vert < 1$, and in this case we will identify explicit characterisations for boundedness and compactness of $\wcomp$ in terms of the order and type of the multiplier $\psi$.

We recall that $\mcomp$ is constant when $\vert \lambda \vert =1$, thus for compact $\wcomp$, it necessarily follows from \eqref{2.2}  that $\vert \lambda \vert< 1$. 
However, it is important to realise that the converse does not hold in general and, 
in contrast to the case of unweighted composition operators (see Theorem~\ref{thmA}), 
the condition $\vert \lambda \vert< 1$  is insufficient in itself to 
guarantee compactness of the weighted composition operator $\wcomp$. 
This is illustrated by the following simple example 
which motivates Theorem~\ref{thm1}. (A similar example can be found in \cite[Remark 3.7 (2)]{TK19}.)

\begin{eg} \label{eg:bddNotCpt}
	
	Given $\lambda$ with $ 0 < \lambda < 1$, set $\beta = 1-\lambda^2$,
	so that $0 < \beta < 1$. 
	We also set 
	\[
	\psi(z)=e^{\frac{\alpha}{2}\beta z^2} \quad \textrm{ and } \quad \varphi(z) = \lambda z
	\]
	for $z \in \C$.
	Since $\beta < 1$, the function $\psi$ is in the Fock space $\F$ for $1\leq p \leq \infty$. We compute
	\begin{align*}
	\mcomp& = \sup \left\lbrace \vert \psi(z) \vert^2\, e^{ \alpha \left( \vert \varphi(z) \vert^2 - \vert z\vert^2 \right) }  :  z\in\C \right\rbrace \\
	& = \sup \left\lbrace \vert \psi(z) \vert^2\, e^{\alpha(\lambda^2 - 1)	\vert z\vert^2} : z\in\C \right\rbrace \\
	& = \sup \left\lbrace M_\psi(r)^2\, e^{-\alpha\beta r^2} \colon r >0 \right\rbrace
	\end{align*}
	where $M_\psi(r) \coloneqq \max \left\lbrace \vert \psi(z)\vert  :  \vert z\vert =r \right\rbrace$
	denotes the maximum modulus function of $\psi$. 
	Since $M_\psi(r)^2 = e^{\alpha\beta r^2}$, we see that 
	$M_\psi(r)^2 e^{-\alpha\beta r^2} = 1$ for each $r$. 
	In particular, $\mcomp < \infty$ and so, by Theorem~\ref{thmB}, 
	the weighted composition operator $\wcomp$ is bounded on the Fock space. 
	However, \eqref{2.2} does not hold in this case and so $\wcomp$ does not act compactly 
	on $\F$.   \qed
\end{eg}

The parameter $\beta$ from Example \ref{eg:bddNotCpt} plays a key role in the sequel. 
So for $0 < \abs{\lambda} < 1$ we set
\begin{equation}  \label{defn:beta}
\beta \coloneqq 1- \abs{\lambda}^2,
\end{equation} 
which in particular gives that $0 < \beta < 1$. 

We now show that the order and type of the multiplier $\psi$ can shed some 
further light on when the operator $\wcomp$ is bounded and when it is compact. 
We recall that the \emph{order} $\rho$ of an entire function $f$ is defined as
\[
\rho \coloneqq \limsup_{r\to\infty} \frac{\log \log M_f(r)}{\log r}
\]
and, if $f$ is of finite order $\rho$, its \emph{type} $\sigma$ is defined to be 
\[
\sigma \coloneqq \limsup_{r\to\infty} \frac{\log M_f(r)}{r^\rho}.
\]
For example, every polynomial has order zero, while the function $f(z) = e^{a z^2}$ has order 2 and type 
$\vert a \vert$. 

B\'en\'eteau et al.~\cite{BCK10} previously investigated the order and type of 
functions in the Fock space. 
They showed that each function $f \in \Fh$ has order at most 2, 
and that if $f$ has order 2 then its type is at most $\alpha/2$. 
Conversely, if $f$ has order less than 2, or if it has order 2 and type 
less than $\alpha/2$ then $f$ is in the Fock space $\Fh$. 

However, membership of the Fock space $\Fh$ cannot be characterised in terms 
of order and type alone.  In \cite{BCK10} examples of functions in the liminal case are provided, 
that is functions of order 2 and type $\alpha/2$ that belong to the Fock space $\Fh$, 
and examples of such functions that do not. It is also shown that if 
$f \in \Fh$  has order two and type $\alpha/2$ then 
$f$ must have infinitely many zeros. 

In the general case, if  $f \in \F$ with $1 \leq p \leq \infty$, then $f$ has order less than or equal to 2, and  if $f$ is of order 2, then $f$ must be of type less than or equal to $\alpha/2$ (cf.~\cite[Theorem 2.12]{Zhu12}).

Similarly, condition \eqref{2.2} can  \emph{almost} be characterised in terms 
of the growth of the multiplier $\psi$. 
Our next result demonstrates that there is a balance 
between the  growth of the symbol $\varphi$ and that of the 
multiplier $\psi$, in that the more contractive $\varphi$ is, the faster $\psi$ is permitted to grow. 

\begin{thm}\label{thm1}
	Let $1 \leq p \leq \infty$, and  $\varphi(z) = a + \lambda z$ with $\vert \lambda \vert < 1$. 
	If $\psi$ has order strictly less than 2, or if $\psi$ has order 2 and type strictly less than 
	$\alpha\beta/2$, then $\wcomp$ is a compact weighted composition operator 
	on $\F$ and on $\Fzero$. 
	
	If $\psi$ has order greater than 2, or if $\psi$ has order 2 and type strictly greater than 
	$\alpha\beta/2$, then $\wcomp$ is neither bounded on $\F$ nor on $\Fzero$. 
\end{thm}

\begin{proof}Suppose that $\psi$ is of order strictly less than 2, that is
	\[
	\rho = \rho(\psi) = \limsup_{r\to\infty} \frac{\log \log M_\psi(r)}{\log r}  = 2-2\varepsilon,
	\]
	where $\varepsilon$ is positive. 
	Then there exists $r_0$ such that 
	\[
	\frac{\log \log M_\psi(r)}{\log r} \leq 2-\varepsilon
	\]
	for $r\geq r_0$ and thus 
	$\vert \psi(z)\vert  \leq \exp \left( r^{2-\varepsilon} \right)$
	for $\vert z \vert  = r \geq r_0$. 
	Then $\psi$ belongs to $\F$, $1\leq p < \infty$, and to $\Fzero$.
	Next notice that
	\[
	\vert \varphi(z) \vert^2 -\vert z\vert^2 =
	\vert a+\lambda z\vert^2 -\vert z\vert^2 
	\leq  -\beta \vert z\vert^2 + \vert a\vert^2 + c_1 \vert z\vert, 
	\]
	where $c_1 = 2 \vert a \vert \vert \lambda\vert$. 
	Thus, for $\vert z\vert =r \geq r_0$, 
	\[
	\vert \psi(z)\vert^2 e^{\alpha \left( \vert \varphi(z)\vert^2 - \vert z \vert^2 \right)}  
	\leq c_2 \exp \left( r^2 \left[  -\alpha\beta + 2r^{-\varepsilon } + \alpha c_1/r \right] \right)
	\]
	where $c_2 = e^{\alpha\vert a\vert^2}$. This shows that \eqref{2.2} holds and so it follows from Theorem~\ref{thmC} that $\wcomp$ is compact. 
	
	Now suppose that $\psi$ has order 2 and type $\sigma$
	with $\sigma < \frac{\alpha}{2}\beta$, 
	say $\sigma = \frac{\alpha}{2}\beta - \varepsilon$
	where $\varepsilon$ is positive. 
	Since $\psi$ has order 2 and type $\sigma$, there exists $r_0$ such that 
	\[
	\frac{\log M_\psi(r)}{r^2} \leq \tfrac{\alpha}{2}\beta- \tfrac{\varepsilon}{2}, 
	\quad r \geq r_0.
	\] 
	Thus $M_\psi(r) \leq e^{r^2(\tfrac{\alpha}{2}\beta- \tfrac{\varepsilon}{2})}$, 
	for $ r \geq r_0$, so that $\psi$ belongs to $\F$, $1\leq p < \infty$, and to $\Fzero$. Moreover, 
	$\vert \psi(z)\vert^2\leq e^{r^2(\alpha\beta - \varepsilon)}$
	for $\vert z\vert = r \geq r_0$ and we have that
	\[
	\vert \psi(z)\vert^2 e^{\alpha \left(\vert \varphi(z)\vert^2 - \vert z \vert^2 \right)} 
	\leq c_2 \exp \left( r^2 \left[ -\alpha\beta + 
	(\alpha\beta - \varepsilon) + \alpha c_1/r \right] \right).
	\]
	So we again see that \eqref{2.2} holds and  $\wcomp$ is compact.

	Next, suppose that $\psi$ has order 2 and type $\sigma$ where 
	$\sigma > \frac{\alpha}{2} \beta$, say $\sigma = \frac{\alpha}{2} \beta + \varepsilon$ 
	where $\varepsilon$ is positive.  
	Then there exists a sequence $r_n \to\infty$ such that 
	\[
	\frac{\log M_\psi(r_n)}{r_n^2} \geq \frac{\alpha}{2} \beta  + \frac{\varepsilon}{2}, 
	\]
	that is $M_\psi(r_n) \geq e^{ r_n^2(\tfrac{\alpha}{2} \beta  + \tfrac{\varepsilon}{2})}$.
	Thus, there exists a sequence $\{ z_n \}$, with $\vert z_n \vert = r_n \to \infty$, 
	such that 
	$\vert \psi(z_n) \vert^2  \geq e^{ r_n^2(\alpha\beta + \varepsilon)}$. 
	Then
	\[
	\vert \psi(z_n)\vert^2  e^{ \alpha \left( \vert \varphi(z_n)\vert^2 - \vert z_n \vert^2 \right)} 
	\geq c_2 \exp\left( (\alpha\beta + \varepsilon) r_n^2 
	- \alpha\beta r_n^2 - \alpha c_1 r_n  \right)
	\]
	which is unbounded as $n\to\infty$.  
	Since $\mcomp$ is not finite it follows that the weighted 
	composition operator $\wcomp$ is not bounded on $\F$ nor on $\Fzero$. 
	
	Finally, let $\psi$ have order strictly greater than 2. It follows from \cite[Theorem 2.12]{Zhu12} that 
	$\psi \notin \F$, $1\leq p< \infty$,
	and $\psi \notin \Fzero$, 
	and thus Theorem \ref{thmB} gives that $\wcomp$ is not bounded on $\F$ nor on $\Fzero$.
\end{proof}

When  $\vert \lambda\vert < 1$, Theorem \ref{thm1} gives that  $\wcomp$ can be bounded and non-compact only if $\psi$ has order 2 and type $\alpha\beta/2$.
However, in this case it is not possible to determine based on this information alone whether or not 
the operator $\wcomp$ is compact, or even bounded. 
Indeed, Example \ref{eg:bddNotCpt} shows that such a choice of $\psi$ and 
$\varphi$ may give rise to an operator $\wcomp$ that is 
bounded on the Fock space but not compact.  

The subtlety of this question is further revealed by considering $\psi(z) = z\,e^{\frac{\alpha}{2}\beta z^2}$ and $\varphi(z) = \lambda z$, which gives that
$\mcomp = \infty$ and so $\wcomp$ is not even bounded on $\F$. 
On the other hand,  if we take $\psi(z) = \big( e^{\frac{\alpha}{2}\beta z^2} -1 \big)/z$ then 
\eqref{2.2} holds and $\wcomp$ acts compactly on the Fock space. 

It is, however, possible to give a full description of the zero-free multipliers $\psi$ associated 
with the symbol $\varphi(z) = a + \lambda z$, for $\vert \lambda \vert < 1$. 
For such multipliers, the operator $\wcomp$ is compact if and only if 
$\psi$ has order less than two, or has order two and type strictly less than 
$\alpha\beta/2$, so that a converse to 
the sufficient compactness condition in Theorem~\ref{thm1} holds 
for zero-free multipliers.

\begin{thm}\label{thm2}
	Let $\varphi(z) = a + \lambda z$ with $\vert \lambda \vert < 1$,
	and suppose that $\psi$ is non-vanishing. 
	Consider the weighted composition operator $\wcomp$ acting on $\F$, 
	$1 \leq p \leq \infty$, or on $\Fzero$.
	
	\begin{enumerate}[label=(\alph*), leftmargin=*, itemsep=1ex]
		
		\item \label{item:a} 
		$\wcomp$ is compact
		if and only if $\psi$ has the form 
		\beq\label{SC1}
		\psi(z) = e^{a_0 + a_1 z + a_2 z^2}
		\eeq
		and $\vert a_2 \vert < \frac{\alpha}{2}\beta$.
		
		\item \label{item:b} 
		$\wcomp$ is bounded but not compact if and only if $\psi$ has the form \eqref{SC1}	
		with $\vert a_2 \vert = \frac{\alpha}{2}\beta$, and  
		either $t := a_1 + \alpha\overline{a}\lambda$ is zero or else 
		both  $t \neq 0$ and  $a_2 = - \frac{\alpha}{2}\beta
		\,\frac{t^2}{\abs{t}^2}$. 
	\end{enumerate}
\end{thm}

\begin{proof}
	Before we begin we record some facts that are required in the proof.
	
	Suppose that $\wcomp$ is bounded.
	It then follows from Theorem \ref{thm1} 
	that $\psi$ has order at most 2. Since $\psi$ is non-vanishing, 
	it follows from the Hadamard Factorisation Theorem  that
	\[
	\psi(z) = e^{g(z)}
	\]
	where $g(z) = a_0 + a_1 z + a_2 z^2$ is a polynomial of degree at most 2 (the order of $\psi$). 
	If $a_2 = 0$, then $\psi$ has order strictly less than 2 and  Theorem~\ref{thm1} gives that $\wcomp$ is compact (in agreement with parts \ref{item:a} and \ref{item:b}). 
	
	Thus we may suppose from now on that $a_2 \neq 0$, in which case 
	$\psi$ has order 2 and non-zero type $\vert a_2\vert$.
	Theorem~\ref{thm1} then gives that $\psi$ has type at most 
	$\frac{\alpha}{2}\beta$, and hence  $\vert a_2 \vert \leq \frac{\alpha}{2}\beta$.
	A computation gives that
	\beq\label{SC21}
	\mzcomp = C\, \exp\left( 2\, \Re \left( t z \right)	+ 2\, \Re (a_2 z^2 ) - \alpha\beta \vert z\vert^2  \right),
	\eeq
	where $C = \exp\left( 2\Re a_0 + \alpha \vert a\vert^2 \right)$ 
	and $t = a_1 + \alpha\overline{a}\lambda$. We
	write $a_2 = \vert a_2\vert e^{-2i\theta_2}$, where $0\leq \theta_2<\pi$,
	and we replace $z$ by $ e^{i\theta_2}w$ on the right-hand side of \eqref{SC21} to give 
	$\mzcomp = C\, e^{2h(w)}$ where
	\beq\label{SC21a}
	h(w) := \Re ( t e^{i\theta_2} w ) + \vert a_2\vert\, \Re ( w^2 ) - \frac{\alpha\beta}{2} \vert w\vert^2.
	\eeq
	
	It follows by Theorem~\ref{thmB} that $\wcomp$ is bounded if and only if 
	$\psi \in \F$, $1\leq p\leq \infty$, respectively $\psi \in \Fzero$, and
	$h$ is bounded on $\C$. 
	By Theorem~\ref{thmC} it follows that $\wcomp$ is compact if and only if 
	$\psi \in \F$, $1\leq p\leq \infty$, respectively $\psi \in \Fzero$, and
	$h(w) \to -\infty$ as $\vert w \vert \to \infty$. 
	Since \eqref{defn:beta} gives that $\beta < 1$, the fact that  $\psi \in \F$, $1\leq p\leq \infty$, 
	respectively $\psi \in \Fzero$, is immediate from Remark~\ref{rmk:QuadraticsInFock} and the condition that $\abs{a_2} \leq \frac{\alpha}{2}\beta < \frac{\alpha}{2}$. Thus, in this case $\mcomp$ is the determining factor
	of boundedness or compactness of $\wcomp$.
	
	\medskip\noindent
	\ref{item:a}
	Suppose that $\wcomp$ is compact. 
	It follows that $\psi$ is of the form \eqref{SC1} with $\abs{a_2} \leq \frac{\alpha}{2}\beta$.
	
	If $\abs{a_2} =\frac{\alpha}{2}\beta$, observe for $w= u$ real that \eqref{SC21a} gives 
	\beq\label{SC21real}
	h(u) = u\, \Re ( t e^{i\theta_2}  ).
	\eeq
	which does not tend to $-\infty$ as $\vert u\vert \to \infty$, irrespective of the value of 
	$\Re \big( t e^{i\theta_2}  \big) $. 
	Consequently, when $\wcomp$ is compact 
	it must hold that $\abs{a_2} < \frac{\alpha}{2}\beta$.
	
	\medskip
	
	Conversely, suppose that $\psi(z) = e^{a_0 + a_1 z + a_2 z^2}$ with 
	$\vert a_2 \vert < \frac{\alpha}{2}\beta$.  It follows from \eqref{SC21a}, for $\abs{w} = r$, that
	\[
	h(w) \leq   \abs{t} r 	- \Big( \frac{\alpha\beta}{2} -\vert a_2\vert \Big) r^2,
	\]
	which gives that $h(w) \to -\infty$ as $r \to \infty$, and thus $\wcomp$ is compact. 
	
	\medskip\noindent
	\ref{item:b}
	Suppose that $\wcomp$ is bounded and non-compact with $\psi$ non-vanishing, 
	so it follows that $\psi$ has the form \eqref{SC1} with $\abs{a_2} \leq\frac{\alpha}{2}\beta$.  
	In fact, $\abs{a_2} =\frac{\alpha}{2}\beta$ in this case, 
	since part \ref{item:a} of this theorem gives that $\wcomp$ is compact when 
	$\abs{a_2} < \frac{\alpha}{2}\beta$. 
	Thus
	\beq\label{SC21b}
	h(w) = \Re ( t e^{i\theta_2} w )	+ 
	\frac{\alpha\beta}{2}\left( \Re (w^2) -  \vert w\vert^2 \right).
	\eeq
	When $w = u$ is real in \eqref{SC21b}, we again obtain \eqref{SC21real} and
	hence the necessary condition  
	$\Re \big( t e^{i\theta_2}  \big) = 0$  for $h$ to be bounded on $\C$. 
	The latter is equivalent to the condition that either $t=0$,  
	or  else both $t\neq 0$ and  $e^{-i\theta_2} = \pm i t/\abs{t}$, that is 
	$a_2 = \frac{\alpha}{2}\beta e^{-2i\theta_2} = - \frac{\alpha}{2}\beta \,t^2/\abs{t}^2$.
	
	\medskip
	
	Conversely we assume that $\psi$ is of the form \eqref{SC1}, 
	with $\vert a_2 \vert = \frac{\alpha}{2}\beta$, 
	and that either $t=0$ or else both  $t \neq 0$ and 
	$a_2 = - \frac{\alpha}{2}\beta	\,t^2/\abs{t}^2$.
	
	Part \ref{item:a} of this theorem gives that $\wcomp$ cannot be compact, 
	so we only need to check that $h$ is bounded on $\C$. 
	
	We first check the case $t = 0$.  By \eqref{SC21b} we have that $2h(w) = \alpha\beta \big(\Re(w^2) - \vert w\vert^2\big)
	\leq 0$, $w \in \C$.
	
	Next we check the case $t \neq 0$ and $a_2 = - \frac{\alpha}{2}\beta	\,t^2/\abs{t}^2$. 
	Then,
	$t e^{i\theta_2} = \pm i \vert t \vert$ is purely imaginary,
	say  $t e^{i\theta_2} = iy$, for some $y \in \R$. 
	Writing $w=u+iv$, \eqref{SC21b} becomes  
	$ 
	h(w) =  -yv -\alpha\beta v^2,
	$
	which shows that $h$ is bounded on $\C$, and hence $\wcomp$ is also bounded in this case.
\end{proof}

\section{Iterates of $\wcomp$ for $\abs{\lambda} < 1$ and $\psi$ zero-free}
\label{sec:orbits}

In this section, we consider the iterates of weighted composition operators $\wcomp$ that are 	bounded and non-compact.  We find different behaviour depending on whether $\lambda$ is real or complex. When $\lambda$ is real, this further depends on whether $p$ is finite or $p=\infty$. 
The latter reflects the fact that the function $e^{\alpha z^2 / 2}$ 
belongs to $\Finfty$, but not to $\F$ for $p<\infty$. 

We begin by deriving an asymptotic formula for the iterates of $\wcomp$. 
We denote  the fixed point of $\varphi$ by
\begin{equation} \label{fixedPoint}
z_0 = a/(1-\lambda)
\end{equation}
and we write the $n$-fold composition of $\varphi$ with itself as
\begin{equation*}
\varphi_n = \underbrace{\varphi \circ \varphi \circ \cdots \circ 
	\varphi}_{n\textrm{-times}},
\end{equation*}
where $\varphi_0$ is the identity.  
Then $\varphi_n$ is given explicitly as 
\beq\label{phik}
\varphi_n(z) = a \left(\frac{1-\lambda^n}{1-\lambda}\right) + \lambda^n z = z_0(1-\lambda^n) + \lambda^n z.
\eeq

\begin{lemma}\label{propiterate}
	Let $\varphi(z) = a + \lambda z$ where $\vert \lambda \vert  < 1$, and let $\psi$ be zero-free and of the form \eqref{SC1}. 
	Suppose that $\wcomp$ is bounded on $\F$, $1\leq p\leq \infty$.
	Then there exist convergent sequences of complex numbers 
	$\{c_{0,n}\}$ and $\{c_{1,n}\}$, with limits $c_0$ and $c_1$ respectively, 
	such that, for $f \in \F$ and $n \geq 1$, 
	\beq\label{S5}
	\left(\wcomp^n f\right)(z) = \psi(z_0)^n \exp\Big( c_{0,n} + c_{1,n} z + 
	a_2 \frac{1-\lambda^{2n}}{1-\lambda^2} z^2\Big) \, f\left(\varphi_n(z) \right),
	\eeq
	where $z_0$ is the fixed point of $\varphi$ as given by \eqref{fixedPoint}.
	Moreover, 
	\beq\label{C1}
	c_1 = \frac{1}{1-\lambda}\Big( a_1 +  \lambda a\, \frac{2a_2}{1-\lambda^2} \Big).
	\eeq
\end{lemma}

\begin{proof}
	We first note that while it is the quadratic term in the exponential of \eqref{S5} that primarily concerns us, the factor $\psi(z_0)^n$ is also important.
	
	We write $g(z) = a_0 + a_1 z + a_2 z^2$ for the quadratic term in \eqref{SC1}.
	Then, for $f\in\F$ and using the fact that $\psi(z_0) = e^{g(z_0)}$, it follows that
	\begin{align}
	\left( \wcomp^n f \right)(z) & = \left( \,\prod_{k=0}^{n-1}\psi\left( \varphi_k(z) \right) \,\right) f\left( \varphi_n(z) \right)\nonumber\\
	& = \left( \,\prod_{k=0}^{n-1}e^{g( \varphi_k(z) )}\,\right) f\left( \varphi_n(z) \right)\nonumber\\
	& = \psi(z_0)^n\, \left( \,\prod_{k=0}^{n-1}e^{g( \varphi_k(z) ) - g( z_0 ) } \, \right)
	f\left( \varphi_n(z) \right).
	\label{wn}
	\end{align} 
	Next we note that 
	\begin{align}
	\prod_{k=0}^{n-1}e^{g( \varphi_k(z) ) - g(z_0) }&  = 
	\prod_{k=0}^{n-1} e^{a_1 (\varphi_k(z) - z_0 ) + 
		a_2 ( \varphi^2_k(z) - z_0^2 )}\nonumber\\
	& =
	\exp\left[a_1 \sum_{k=0}^{n-1} \left(\varphi_k(z) - z_0 \right) + 
	a_2 \sum_{k=0}^{n-1} \left( \varphi^2_k(z) - z_0^2 \right)\right].
	\label{product}
	\end{align}
	Now it follows from \eqref{phik} that $\varphi_k(z) - z_0 = \lambda^k(z-z_0)$ and
	\[
	\varphi^2_k(z) - z_0^2 = \left( \varphi_k(z) - z_0 \right)  \left( \varphi_k(z) + z_0 \right)
	= \lambda^{2k}(z-z_0)^2 + 2\lambda^k z_0(z-z_0).
	\]
	Thus,
	\begin{align} 
	\sum_{k=0}^{n-1} \left( \varphi_k(z) - z_0 \right)
	&= \frac{1-\lambda^n}{1-\lambda} (z-z_0) \label{sumphi}
	\shortintertext{and}
	\sum_{k=0}^{n-1}\left(  \varphi^2_k(z) - z_0^2 \right) 
	&= \frac{1-\lambda^{2n}}{1-\lambda^2} (z-z_0)^2 + 2\frac{1-\lambda^n}{1-\lambda}z_0 (z-z_0). \label{sumphi2}
	\end{align}
	Formula \eqref{S5} then follows from \eqref{wn}, \eqref{product}, \eqref{sumphi} and \eqref{sumphi2}, 
	with explicit expressions for the coefficients $c_{0,n}$ and $c_{1,n}$. 
	The only quadratic term comes from \eqref{sumphi2}. 
	Computing the term $c_{1,n}$, taking the limit as $n \to \infty$ (so that 
	$\lambda^n \to 0$), and finally replacing $z_0$ by $a/(1-\lambda)$ leads to \eqref{C1}. 
\end{proof}

The simplest case when studying the iterates of $\wcomp$ is when $\lambda$ is not real. 

\begin{thm}\label{lambdacomplex}
	Let $\varphi(z) = a + \lambda z$ where $\vert \lambda \vert < 1$ and $\lambda$ is not real,
	and let $\psi$ be zero-free and of the form \eqref{SC1}. 
	Suppose that $\wcomp$ is bounded on $\F$, $1\leq p \leq \infty$.  
	Then $\wcomp^2$ is a compact weighted composition operator on $\F$. 
\end{thm}

\begin{proof}
	It follows from \eqref{wn}, \eqref{product}, \eqref{sumphi} and \eqref{sumphi2} that, for any $f \in \F$,  
	\[
	\wcomp^2f(z) = \psi(z_0)^2 \, \exp\left[ P(z) \right]\, f\left( a(1+\lambda) + \lambda^2 z \right)
	\]
	where, explicitly,
	\[
	P(z) = a_1(1+\lambda) (z-z_0) + a_2 \left[ (1+\lambda^2)(z-z_0)^2 + 2(1+\lambda) z_0 (z-z_0) \right].
	\]
	This implies that $\wcomp^2$ is itself a weighted composition operator 
	with symbol $\widetilde\varphi(z) = a(1+\lambda) + \lambda^2z$ 
	and multiplier $\widetilde\psi(z) = \psi(z_0)^2 \, \exp\left[ P(z) \right]$. 
	Note that $\tilde{\beta} = 1-\vert \lambda^2\vert^2 = 1- \vert \lambda\vert^4$.
	
	Since $\vert a_2\vert \leq \alpha (1-\vert \lambda \vert^2)/2$, the modulus of the coefficient $\tilde a_2$ of $z^2$ in the quadratic $P$ satisfies
	\[
	\abs{ \tilde a_2 } = \abs{ a_2 } \, \abs{ 1+\lambda^2}  
	\leq  \tfrac{\alpha}{2}\big( 1- \abs{\lambda}^2 \big) \, \abs{ 1+\lambda^2}
	< \tfrac{\alpha}{2} \big( 1- \abs{\lambda}^4 \big) = \tfrac{\alpha}{2} \tilde\beta. 
	\]
	The last strict inequality is a consequence of the assumption that $\lambda$ is not real.
	So it now follows from Theorem~\ref{thm2} that $\wcomp^2$ is compact.
\end{proof}

We now consider the case when $\lambda$ is real. As mentioned earlier, 
the behaviour of the iterates of $\wcomp$ depends on $p$. We first consider the case of finite $p$.

\begin{thm}\label{propreal}
	Let $\varphi(z) = a + \lambda z$, where $-1 < \lambda < 1$ and $\lambda \neq 0$.
	Let $\psi$ be zero-free and of the form \eqref{SC1} with $\vert a_2 \vert= \tfrac{\alpha}{2}\beta$. 
	Suppose that $\wcomp$ is bounded on $\F$, $1\leq p < \infty$.
	If $f \in \F$ with $f(z_0) \neq 0$, where $z_0$ is the fixed point of $\varphi$, then
	\beq\label{S7}
	\vert \psi(z_0)\vert^{-n}\Vert \wcomp^n f \Vert \longrightarrow  \infty, \textrm{ as } n \to \infty.
	\eeq
\end{thm}

\begin{proof}
	Since $f(z_0) \neq 0$, we may choose $r_0>0$ so that 
	$\vert f(z) \vert \geq \frac{1}{2} \vert f(z_0) \vert =c$ 
	for $z\in D(z_0,r_0)$, the disc of radius $r_0$ centred at $z_0$. Note that $c$ is positive. 
	Set $c_{2,n} = a_2 (1-\lambda^{2n})/(1-\lambda^2)$ and 
	$a_2 = \tfrac{\alpha}{2}\beta e^{-2i\theta_2}$, with $0\leq \theta_2<\pi$, 
	so that 
	\[
	c_{2,n} = \frac{\alpha}{2} (1-\lambda^{2n}) e^{-2i\theta_2}. 
	\]
	Here we used the fact that $\lambda$ is real, which gives that $\beta = 1-\vert \lambda\vert^2 = 1-\lambda^2$. 
	
	We suppose, as we may  by Lemma~\ref{propiterate}, 
	that $\abs{ c_{0,n}} \leq 2\abs{c_0}$  for $n$ large enough. 
	Then, by \eqref{S5}, 
	\begin{align}
	\vert \psi(z_0)\vert^{-n} \Vert \wcomp^n f \Vert& = 
	\left\Vert \exp\left( c_{0,n} + c_{1,n} z + c_{2,n} z^2\right)\, 
	f\left(\varphi_n(z)\right)\right\Vert\nonumber\\
	& \geq  e^{-2 \abs{c_0}}
	\left\Vert \exp\left( c_{1,n} z + 	c_{2,n}  z^2 \right)	f\left(\varphi_n(z)\right)\right\Vert.
	\label{first}
	\end{align}
	Next, for $1\leq p < \infty$, 
	\begin{align}
	&\left\lVert \exp\left( c_{1,n} z + c_{2,n} z^2\right) f\left(\varphi_n(z)\right) \right\rVert^p \nonumber \\
	& = \frac{\alpha p}{2\pi} \int_\C \exp\left[ p\,\Re(c_{1,n} z)  + p\,\Re(c_{2,n} z^2)\right] \left\lvert f\left( \varphi_n(z)\right) \right\rvert^p	e^{-\tfrac{\alpha p}{2} \vert z \vert^2}\dm(z)\nonumber\\
	& \geq \frac{c^p \alpha p}{2 \pi} \int_{D_n} \exp\left[ p\,\Re(c_{1,n} z)  + p\,\Re(c_{2,n} z^2)\right]\,
	e^{-\tfrac{\alpha p}{2} \vert z\vert^2}\dm(z)\label{S8a}\\
	& = \frac{c^p \alpha p}{2 \pi} \int_{D_n} 
	\exp\left[ p\,\Re(c_{1,n} e^{i\theta_2}w)  + p\, \vert c_{2,n} \vert\,\Re(w^2)\right]\,
	e^{-\tfrac{\alpha p}{2} \vert w\vert^2}\dm(w) \label{S8}
	\end{align}
	for sufficiently large $n$.  The last step follows by substituting $z = e^{i\theta_2}w$, and we let $D_n = D(0, \tfrac{1}{2} r_0 \vert \lambda \vert^{-n})$ denote the disc of radius $\tfrac{1}{2} r_0 \vert \lambda \vert^{-n}$ centred at the origin. 
	To see that inequality \eqref{S8a} holds, let $z\in D_n$.  It follows by \eqref{phik} that
	\[
	\left\vert \varphi_n(z) - z_0 \right\vert 
	=  \vert\lambda\vert^n \vert z-z_0\vert
	\leq \vert\lambda\vert^n \vert z\vert +\vert\lambda\vert^n \vert z_0\vert
	\leq \tfrac{1}{2} r_0 + \vert \lambda \vert^n \vert z_0\vert < r_0
	\]
	for sufficiently large $n$. 
	Thus, $\varphi_n(z) \in D(z_0,r_0)$ for $z \in D_n$, and hence 
	$\abs{ f(\varphi_n(z))} \geq c$. 
	
	Notice that
	\begin{align*}
	\exp & \left[ p\,\Re(c_{1,n}  e^{i\theta_2}w)  + p \abs{c_{2,n}} \Re(w^2) \right]  \nonumber \\
	&\qquad + \exp\left[ p\,\Re(c_{1,n} e^{i\theta_2}(-w))  + p \abs{c_{2,n}} \Re((-w)^2) \right] \nonumber\\
	& = \left(\exp\left[ p\,\Re(c_{1,n} e^{i\theta_2}w) \right] 
	+ \exp\left[ - p\,\Re(c_{1,n} e^{i\theta_2}w)  \right] \right)\, \exp \left[ p\, \vert c_{2,n} \vert\,\Re(w^2) \right]\nonumber\\
	& \geq 2 \exp \left( p \abs{c_{2,n}} \Re(w^2)\right).	
	\end{align*}
	We set
	$S \coloneqq \left\lbrace  w = re^{it} \,  : \,  - \tfrac{\pi}{4} < t < \tfrac{\pi}{4},\, r > 0 \right\rbrace$.
	Summing the integrals over $D_n\cap S$ and $D_n \cap (-S)$, we find that 
	\begin{align}
	\int_{D_n} &
	\exp\left[ p\,\Re(c_{1,n} e^{i\theta_2}w)  + p\, \vert c_{2,n} \vert\,\Re(w^2)\right]\,
	e^{-\tfrac{\alpha p}{2} \vert w\vert^2}\dm(w)\nonumber\\
	& \geq 
	2\int_{D_n\cap S} \exp\left[ p\, \vert c_{2,n} \vert\,\Re(w^2) -\tfrac{\alpha p}{2}  \vert w\vert^2\right]\dm(w)\nonumber\\
	& = 2 \int_{D_n\cap S} \exp\left[ \tfrac{\alpha p}{2}(1-\lambda^{2n})\,\Re(w^2) 
	-\tfrac{\alpha p}{2}\vert w\vert^2\right]\dm(w).
	\label{DnS}
	\end{align}
	Set
	\begin{align}
	g_n(w) &= \exp\left[ \tfrac{\alpha p}{2}(1-\lambda^{2n})\,\Re(w^2) -
	\tfrac{\alpha p}{2} \vert w \vert^2\right]\, \mathbbm{1}_{D_n\cap S}
	\shortintertext{and}
	g(w) &= \exp\left[ \tfrac{\alpha p}{2}\, \Re(w^2) - 
	\tfrac{\alpha p}{2} \vert w \vert^2\right]\, \mathbbm{1}_S. 
	\end{align}
	Since $\Re(w^2) > 0$ on the sector $S$ and  $1-\lambda^{2n}$ 
	increases to 1, 
	the functions $g_n$ increase monotonically to $g$ on $S$. 
	Hence, by monotone convergence, 
	\begin{align*}
	\int_{D_n\cap S} &  \exp\Big[ \tfrac{\alpha p}{2}(1-\lambda^{2n})\, \Re(w^2) 
	- \tfrac{\alpha p}{2}\vert w \vert^2 \Big]\dm(w)\\
	&\xrightarrow{n\to\infty} 
	\int_{ S} \exp\left[ \tfrac{\alpha p}{2} \Re(w^2) 
	- \tfrac{\alpha p}{2} \vert w\vert^2\right]\dm(w)
	= \int_S \exp \left[ - \alpha p v^2 \right] \du\dv.
	\end{align*}
	At the last step we set $w = u+iv$. Since this last integral is infinite, 
	\eqref{first}, \eqref{S8} and \eqref{DnS} together establish \eqref{S7}. 
\end{proof}

The proof of Theorem~\ref{propreal} shows (using the 
notation of Lemma~\ref{propiterate} and recalling that 
$\beta = 1-\lambda^2$ when $\lambda$ is real) that
\[
\Big\Vert \exp\Big(c_{1,n} z + a_2\beta^{-1} (1-\lambda^{2n}) z^2\Big) \Big\Vert \to \infty 
\mbox{ as } n \to \infty.
\]
In contrast, we will see in the case $\lambda$ real and $p=\infty$, that the exponential term in \eqref{S5} in fact converges for the iterates of $\wcomp$.

\begin{thm}\label{prop:realinfinity}
	
	Let $\varphi(z) = a + \lambda z$ where $-1 < \lambda < 1$ and $\lambda \neq 0$. Let $\psi$ be zero-free and of the form \eqref{SC1} where $\vert a_2 \vert= \tfrac{\alpha}{2}\beta$.
	Suppose that $\wcomp$ is bounded on $\Finfty$. 
	Then the function 
	\beq\label{Wlimit}
	F(z) = \exp\left( c_1z + a_2\beta^{-1} z^2 \right)
	\eeq
	belongs to $\Finfty \setminus \Fzero$, where $c_1$ is as in the statement of Lemma~\ref{propiterate}.
	Moreover, for $f\in \Finfty$, $\psi(z_0)^{-n} \wcomp^nf \to Wf$ locally uniformly on $\C$ 
	where 
	\beq\label{Wlimit2}
	(Wf)(z) = c F(z) f(z_0), \quad f \in \Finfty,
	\eeq
	and $c = e^{c_0}$ where $c_0 = \lim_{n\to\infty} c_{0,n}$ in Lemma~\ref{propiterate}. 
\end{thm}

\begin{proof}
	We again write $a_2 = \tfrac{\alpha}{2}\beta e^{-2i\theta_2}$
	with $0\leq \theta_2<\pi$. 
	Since $\lambda$ is real we have that 
	$2a_2 / (1-\lambda^2) = 2a_2\beta^{-1} = \alpha e^{-2i\theta_2}$. 
	
	Let $t = a_1 + \alpha \overline{a} \lambda$.
	By Theorem~\ref{thm2} \ref{item:b}, the assumption that $\wcomp$ is bounded gives that either $t= 0$, 
	or both $t\neq0$ and $a_2 = -\tfrac{\alpha}{2}\beta t^2 / \abs{t}^2$. 
	
	We first consider the case $t=0$, so that $a_1 = -\alpha \overline{a}\lambda$. 
	It follows from \eqref{C1} that 
	\beq\label{C1infinity1}
	c_1 	= \frac{-\alpha \overline{a} \lambda + \alpha a \lambda e^{-2i\theta_2}}{1-\lambda}
	= e^{-i\theta_2}\frac{\alpha\lambda}{1-\lambda} 2i\, \Im\left( a e^{-i\theta_2} \right). 
	\eeq
	By definition of the norm on $\Finfty$,
	\begin{align*}
	\norm{F} & = \sup_{z\in\C} \left\{ 
	\abs{ \exp\left( c_1z + a_2\beta^{-1} z^2 \right) } \, 
	e^{-\tfrac{\alpha}{2} \abs{ z}^2} \right\}\nonumber\\
	& = \sup_{z\in \C} \left\{ \exp\left[ \Re(c_1 z)  + \Re(\tfrac{\alpha}{2} e^{-2i\theta_2} z^2) -\tfrac{\alpha}{2} \abs{ z}^2 \right] \right\}\nonumber\\
	& = \sup_{w\in \C} \left\{ \exp\left[ \Re(c_1 e^{i\theta_2}w)  + 
	\tfrac{\alpha}{2} \left(\Re( w^2 )  -	\abs{w}^2 \right) \right] \right\},\label{normF1}
	\end{align*}
	where we substituted $w = e^{-i\theta_2}z$. 
	By \eqref{C1infinity1}
	\[
	\Re(c_1 e^{i\theta_2}w) = - \frac{2\alpha\lambda}{1-\lambda} \, 
	\Im( a e^{-i\theta_2} )\, \Im(w) = A \, \Im(w),
	\]
	where $A$ is real. Writing $w=u+iv$,
	\[
	A \, \Im(w) + \tfrac{\alpha}{2} \left( \Re( w^2 )  - \abs{w}^2 \right) = Av - \alpha v^2
	\]
	and thus $e^{Av - \alpha v^2}$ is uniformly bounded on $\C$, thereby showing that $\Vert F\Vert < \infty$, as required. 
	However, $\abs{ F(z)} e^{-\tfrac{\alpha}{2} \abs{z}^2} = e^{Av - \alpha v^2} \not\to 0$ as 
	$\abs{z} \to \infty$, so $F \not\in \Fzero$.
	
	\medskip
	
	Next we consider the case $t\neq0$, with $a_2 = -\tfrac{\alpha}{2}\beta t^2 / \abs{t}^2$. 
	The expression \eqref{C1} for $c_1$ becomes, in terms of $t$, 
	\beq\label{C1infinity2}
	c_1 = \frac{1}{1-\lambda} \left( a_1 - \alpha \lambda a \frac{t^2}{\abs{t}^2} \right).
	\eeq
	Then
	\begin{align*}
	\Vert F \Vert 
	& = \sup_{z\in \C} \left\{ \exp\left[ \Re(c_1 z)  - 
	\Re(\tfrac{\alpha}{2}\tfrac{t^2}{\abs{t}^2}z^2) -
	\tfrac{\alpha}{2} \vert z\vert^2 \right] \right\}\nonumber\\
	& = \sup_{w\in \C} \left\{ \exp\left[ \Re(c_1 \tfrac{\abs{t}}{t} w)  - 
	\tfrac{\alpha}{2} \left[ \Re( w^2 )  +
	\abs{w}^2 \right] \right] \right\},\label{normF2}
	\end{align*}
	where we substituted $w = t z /\abs{t}$. 
	Starting from \eqref{C1infinity2} and replacing $a_1$ by 
	$t- \alpha \overline{a}\lambda$, we find that 
	\begin{align*}
	c_1 \frac{\abs{t}}{t} & = 
	\frac{1}{1-\lambda} \left[ t- \alpha \overline{a}\lambda - 
	\alpha \lambda a \frac{t^2}{\abs{t}^2} \right]\,\frac{\abs{t}}{t}
	\\
	& = \frac{1}{\abs{t} (1-\lambda)} 
	\left[ \abs{t}^2 - 2\alpha\lambda\,\Re(at) \right] = B, 
	\end{align*}
	where $B$ is real. 
	Writing $w = u + iv$, we find that
	\begin{equation*}
	B \, \Re(w) - \tfrac{\alpha}{2} \left( \Re( w^2 )  + \abs{w}^2 \right) = Bu - \alpha u^2
	\end{equation*}
	and thus $e^{Bu - \alpha u^2}$ is uniformly bounded on $\C$, and thus $\norm{F}$ is finite in this case. 
	However, once again we have that 
	$\abs{ F(z)} e^{-\tfrac{\alpha}{2} \abs{z}^2} = e^{Bu - \alpha u^2}  \not\to 0$ 
	as $\abs{z} \to \infty$ and thus $F \not\in \Fzero$.

	Since $\Vert F \Vert_{\Finfty} < \infty$, the weighted point evaluation $W$ given by \eqref{Wlimit2}
	is bounded on $\Finfty$. 
	Comparing the expression \eqref{S5} for $\wcomp^n f$ and the definition of the operator $W$ 
	involving the function $F$ of \eqref{Wlimit}, and noting that $\varphi_n(z) \to z_0$ locally uniformly, 
	it follows that $\psi(z_0)^{-n} \wcomp^nf \to Wf$ locally uniformly on $\C$, for $f \in \Finfty$.
\end{proof}

\section{Non-supercyclic weighted composition operators on Fock spaces} \label{sec:SuperCyFock}

As an application of the results from the previous sections, we now provide a proof that 
weighted composition operators acting on Fock spaces cannot be supercyclic.
In contrast, we recall the elegant result of Bourdon and Shapiro~\cite{BS90,BS97} who proved that a composition operator, induced by an automorphism $\varphi$ of the disc, is hypercyclic on the Hardy space if and only if $\varphi$ has no fixed points in the disc (cf.~\cite[Theorem 4.48]{GEP11}). 
Thus weighted composition operators may have contrasting dynamical behaviours depending on the setting.

The central notion of linear dynamics is hypercyclicity, and we recall that a bounded linear operator $T$ acting on an infinite-dimensional and separable Banach space $X$ is \emph{hypercyclic} 
if there exists $x \in X$ with a dense $T$-orbit in $X$.  The operator $T$ is said to be \emph{supercyclic} if there exists  $x \in X$, 
called a supercyclic vector for $T$, such that its projective orbit 
\begin{equation*}
\{ \zeta \, T^n x : \zeta \in \C, n \in \N \}
\end{equation*}
is dense in $X$.  
It is well known that the class of hypercyclic operators is strictly contained in the class of supercyclic operators. 
We remark that since $\Finfty$ is non-separable, we will only consider weighted composition operators acting on $\F$, for $1 \leq p < \infty$, and $\Fzero$.

In the Fock space setting,  linear dynamical properties of composition operators have been investigated, for instance, in \cite{GI08} and \cite{SMB20}.  In particular, it was proven by Jiang et al.~\cite{JPZ17} and Mengestie and Seyoum~\cite{MS20} that composition operators acting on Fock spaces cannot be supercyclic. Our next result proves that this is also the case for weighted composition operators.
\begin{thm}\label{thm3}
	Let $\wcomp$ be a bounded weighted composition operator acting on the Fock space  $\Fzero$, or on $\F$, for $1 \leq p < \infty$.
	Then $\wcomp$ cannot be supercyclic.
\end{thm}

For a  bounded weighted composition operator $\wcomp$, 
with $\varphi(z) = a + \lambda z$ for $\vert \lambda \vert  \leq 1$, 
verification of Theorem~\ref{thm3} comes down to checking a number of cases
depending on the value of $\lambda$. 
For $\abs{\lambda} = 1$ we will see that $\wcomp$ is a constant multiple of an isometry and hence cannot be supercyclic (the proof of this fact is postponed until Section \ref{sec:paranormal}).
When $\abs{\lambda} < 1$, the behaviour of $\wcomp$ further depends on the form of the multiplier $\psi$, on whether $\lambda$ is real or complex, and on the setting $\F$, $1 \leq p < \infty$, or $\Fzero$.
This is unavoidable since, as we have already seen, the orbits of $\wcomp$ have this dependence.

An initial observation is that $\wcomp$ can only be supercyclic if the multiplier $\psi$ is zero-free.  To see this assume that $\psi$ has a zero at $w$, say, and observe for any $f$ and for $n\geq 1$ that
\[
\left(\wcomp^n f \right)(w) =  \psi(w)  \wcomp^{n-1}(f\circ \varphi)(w) =0
\]
and hence the projective orbit of any $f$ is contained in $\spn{f} \cup {\rm ker}\,(k_{w}z)$,  which  cannot be projectively dense in $\F$, for $1 \leq p < \infty$, or in $\Fzero$,
and hence $\wcomp$ cannot be supercyclic. 
We may therefore assume that the multiplier $\psi$ has the form \eqref{SC1}.

Before we present the proof of Theorem~\ref{thm3}, 
we recall some well known necessary conditions for supercyclicity.

The existence of eigenvalues plays a significant role in the theory of linear dynamics. 
A standard result states if an operator $T$ 
acting on a separable Banach space is hypercyclic, or if $T$ is compact and supercyclic, 
then the point spectrum of the adjoint $T^*$ is empty 
(cf.~\cite[p.~29 and Propositions~1.17 and 1.26]{BM09}).   

Another useful necessary condition for supercyclicity is the following corollary of the Angle Criterion, which can be found in \cite[Corollary~9.3]{BM09}.  Let $T$ be an operator acting on a separable Banach space $X$.
Assume for all $x$ in some non-empty open subset $U \subset X$ that one can find 
a non-zero linear functional $x^* \in X^*$ such that
\beq\label{S1}
\lim_{n\to\infty}\frac{\Vert (T^*)^n x^*\Vert}{\Vert T^n x \Vert} =0.
\eeq
Then $T$ is not supercyclic. 

We are now ready to prove that the Fock spaces do not support supercyclic weighted composition operators. The fact that no such operator is hypercyclic is straightforward, as can be seen in the proof of Theorem \ref{thm3}. 

\begin{proof}[Proof of Theorem~\ref{thm3}] 
	Suppose that $\wcomp$ is a  bounded weighted composition operator acting on the Fock space $\F$, for $1\leq p < \infty$, or on $\Fzero$.  
	Then $\varphi(z) = a + \lambda z$ and  $\vert \lambda \vert  \leq 1$.  
	We also assume, as we may, the necessary condition that $\psi$ is non-vanishing.
	
	We first consider the case $\abs{\lambda} = 1$. 
	Theorem \ref{thm:ConstIsometry} gives that  in this case $\wcomp$ is a constant multiple of an isometry, so it follows from results by Bourdon~\cite{Bou97} (cf.~\cite[p.~159]{GEP11}) that $\wcomp$ cannot be supercyclic (and therefore not hypercyclic). 
	
	\medskip
	
	We next assume that $\abs{\lambda} < 1$.  
	Then the function $\varphi$ has a fixed point at $z_0 = a/(1-\lambda)$ and 
	\eqref{wstar} gives that $\psi(z_0)$ is an eigenvalue of the adjoint $\wcomp^*$, 
	with corresponding eigenvector $k_{z_0}$.  
	(At this point we can conclude that $\wcomp$ cannot be hypercyclic.)
	
	Since $\wcomp^*$ possesses an eigenvalue, it follows from the above mentioned spectral condition that compact $\wcomp$ cannot be supercyclic (we recall that this includes the case $\lambda = 0$ by Remark \ref{rmk:lambdaEqZero}). 
	
	It thus follows from Theorem~\ref{thm2} that we are left to investigate the case when $0 < \abs{\lambda} < 1$, 
	for $\psi$ of the form \eqref{SC1} with $\vert a_2 \vert = \tfrac{\alpha}{2}\beta$ 
	(that is, $\wcomp$ is a bounded and non-compact weighted composition operator). 
	
	In the case $\lambda$ is real and $p$ is finite, we employ the above mentioned corollary of the Angle Criterion.  We choose as our linear functional the point evaluation $k_{z_0}$, where $z_0$ denotes the fixed point of $\varphi$. 
	Iterating \eqref{wstar} gives that 
	\beq\label{a21}
	\norm{ \left( \wcomp^* \right)^nk_{z_0}}
	= \norm{ \psi(z_0)^n \, k_{z_0}}
	= \abs{ \psi(z_0)}^n \, e^{\alpha \abs{ z_0}^2/2}.
	\eeq
	It follows from \eqref{phik} that $\varphi_k(z) \to z_0$ as $k \to \infty$ for each $z$, so we set 
	\beq\label{S6}
	U = \left\{ f \in \F  \, : \,  f(z_0) \neq 0\right\}.
	\eeq
	Note that $U = \delta_{z_0}^{-1}\left( \C\setminus \{0\} \right)$, where $\delta_{z_0}$ is the point evaluation at $z_0$. 
	The point evaluation is a continuous linear functional and $\C\setminus\{0\}$ 
	is an open set, so $U$ is an open subset of $\F$.
	
	Observe that it follows from \eqref{a21} and Theorem \ref{propreal}, for any $f \in U$, that 
	\begin{equation*}
	\frac{\norm{ ( \wcomp^* )^n k_{z_0}}}{\norm{\wcomp^n f}} = \frac{\abs{ \psi(z_0)}^n \, e^{\alpha \abs{ z_0}^2/2}}{\norm{\wcomp^n f}} \longrightarrow 0, 	\textrm{ as } n \to \infty.
	\end{equation*}
	Thus \eqref{S1} holds and it follows from the corollary of the Angle Criterion that $\wcomp$ cannot be supercyclic in this case.
	
	\medskip
	
	Next we consider the case of the space $\Fzero$ and when $\lambda$ is real. 
	We cannot apply the Angle Criterion in this case since the norm of $F$ is finite, so a separate argument is required.
	Curiously, this case seems to require the most work.
	Suppose, to the contrary, that $f_0 \in \Fzero$ is a supercyclic vector for 
	$\wcomp$ on $\Fzero$.  
	It is then also a supercyclic vector for $\twcomp = \psi(z_0)^{-1} \wcomp$ on $\Fzero$.
	We note that a necessary condition on $f_0$ is that $f_0(z_0) \neq 0$. If not, 
	$( \twcomp^n f_0 )(z_0) = 0$ for $n \geq 1$ so that the projective orbit of 
	$f_0$ under $\twcomp$ is contained in $\spn{f_0} \cup \ker(k_{z_0})$. 
	
	It follows from Theorem~\ref{prop:realinfinity} that $\twcomp^n f \to W f$ 
	in the topology of local uniform convergence on $\C$ for each $f\in \Finfty$.
	Recall from \eqref{Finftynorm} that the norm on $\Finfty$ is given by a supremum, so we may choose $z=z_0$ in this supremum and conclude that 
	\[
	\big\Vert \twcomp^n f_0 \big\Vert \geq 
	\abs{  \exp\left[  c_{0,n} + c_{1,n} z_0 + 
		a_2 \beta^{-1}(1-\lambda^{2n})z_0^2 \right] } \, \abs{ f_0(z_0)}
	\, \exp(-\tfrac{\alpha}{2} \abs{ z_0}^2) 
	\]
	which is in turn greater, for all $n$, than $c$, for a suitable positive constant $c$.
	
	Suppose that there exist sequences $\{n_k\}_k \subset \N$ and 
	$\{ \zeta_{n_k} \}_k \subset \C$  such that 
	\beq\label{last}
	\zeta_{n_k}\, \twcomp^{n_k} f_0 \to g \ \mbox{ in } \Fzero,
	\eeq
	so that $g \in \Fzero$ lies in the closure of the projective orbit of $f_0$ 
	under $\twcomp$.
	We claim that $g$ must be 0 and, if so, $f_0$ cannot be a supercyclic vector for 
	$\twcomp$ on $\Fzero$. 
	Suppose, to the contrary, that $g \neq 0$. 
	A particular consequence of \eqref{last} is that 
	\[
	\norm{g}  = \lim_{k\to\infty} \big\Vert\zeta_{n_k}\, \twcomp^{n_k} f_0 \big\Vert = 
	\lim_{k\to\infty} \abs{ \zeta_{n_k}} \, \big\Vert \twcomp^{n_k} f_0 \big\Vert. 
	\]
	Since $\big\Vert\twcomp^{n_k} f_0 \big\Vert  \geq c$, it follows for all sufficiently large $k$ that
	$\abs{ \zeta_{n_k}} \leq 2 \norm{g}/c \, = C_0$. 
	Convergence in $\Finfty$ implies locally uniform convergence, thus for each positive $R$
	\beq\label{lastone}
	\sup_{\abs{z} \leq R} \abs{ \zeta_{n_k}\, (\twcomp^{n_k} f_0)(z) - g(z) }
	\to 0 
	\eeq
	as $k \to \infty$.  Consequently
	\begin{align*}
	\sup_{\abs{z} \leq R} & \abs{ \zeta_{n_k} (W f_0)(z) - g(z)}\\
	& \leq \sup_{\abs{ z} \leq R} \abs{ \zeta_{n_k}}\, \abs{(W f_0)(z) - (\twcomp^{n_k}f_0)(z)} +
	\sup_{\abs{ z } \leq R}  \abs{ \zeta_{n_k}\, (\twcomp^{n_k} f_0)(z) - g(z)}\\
	& \leq C_0\, \sup_{\abs{ z } \leq R} \abs{ (W f_0)(z) - (\twcomp^{n_k}f_0)(z)}+
	\sup_{\abs{ z } \leq R}  \abs{ \zeta_{n_k}\, (\twcomp^{n_k} f_0)(z) - g(z) }.
	\end{align*}
	By the local uniform convergence of $\twcomp^n f_0$ to $W f_0$ and by 
	\eqref{lastone} we see that
	$ \zeta_{n_k}\, (W f_0)(z)  \to g(z)$
	locally uniformly, so that $g \in \spn{W f_0}$.
	But $g \in \Fzero$ and $W f_0 \in \Finfty\setminus \Fzero$, 
	implying that $g=0$ and contradicting the assumption that 
	$g\neq 0$. 
	
	\medskip
	
	The final case to check is when $\lambda$ is not real. 
	We first recall a result of Ansari~\cite[Theorem 2]{Ans95} which states that any power of a supercyclic operator is itself supercyclic with the same supercyclic vectors.
	Next we note that Theorem~\ref{lambdacomplex} gives that $\wcomp^2$ is a compact weighted composition operator on $\F$, for $1\leq p < \infty$, 
	and on $\Fzero$.  From the above discussion, we know that $\wcomp^2$ cannot be supercyclic and it thus follows that $\wcomp$ itself cannot be supercyclic in this case.   
\end{proof}

\subsection{The operator $\wcomp$ when  $\vert\lambda\vert = 1$}
\label{sec:paranormal}

When $\wcomp$ is bounded on $\F$, $1 \leq p \leq \infty$, we know from Theorem \ref{thmB} and \eqref{defn:formOfPhi} that the symbol $\varphi$ is of the form $\varphi(z) = a + \lambda z$, with $\vert \lambda \vert \leq 1$.
For $\abs{ \lambda} = 1$ it follows from \eqref{defn:formOfPsi} that the multiplier $\psi$ must be of the form $\psi(z) = \psi(0) \exp(-\alpha \overline{a} \lambda z)$, while the range of potential multipliers $\psi$ is much wider when $\vert \lambda \vert < 1$.
In this section we consider the case $\abs{ \lambda} = 1$ and we will show that $\wcomp$ is a constant multiple of an isometry and cannot be supercyclic. We also make some observations on the paranormality of $\wcomp$ in this case.

The following theorem uses the translation operator $T_b \colon \F \to \F$, for $1\leq p \leq \infty$, 
which is defined for any $b \in \C$ as
\begin{equation*}
(T_b f)(z) =  e^{-\alpha \overline{b} z - \frac{\alpha}{2}\abs{b}^2} f(b+z),
\end{equation*}
where  $z \in \C$ and $f \in \F$.  
We note that $T_b$ is an isometry on the Fock spaces $\F$ (cf.~\cite[Proposition 2.38]{Zhu12}).
For $\abs{\lambda} = 1$, we compose $T_b$ with a rotation by $\lambda$ 
to define the operator $T_{\lambda, b} \colon \F \to \F$ by
\begin{equation} \label{defn:Tisometry}
(T_{\lambda, b} f)(z) =  e^{-\alpha \overline{b} \lambda z - \frac{\alpha}{2}\abs{b}^2} f(b + \lambda z).
\end{equation}
It can easily be seen that $T_{\lambda, b}$ is also an isometry.

We note that an isometry cannot be supercyclic, which follows in the Hilbert space setting by Hilden and Wallen~\cite{HW73}, and in the Banach space setting by Ansari and Bourdon~\cite{AB97}.

\begin{thm}  \label{thm:ConstIsometry}
	Let $\wcomp$ be a bounded weighted composition operator, with $\varphi(z) = a + \lambda z$ and $\abs{\lambda} = 1$, on the Fock space $\F$, $1 \leq p \leq \infty$, or on $\Fzero$.  Then $\wcomp$ is a constant multiple of an isometry and cannot be supercyclic.
\end{thm}

\begin{proof}
	For $\abs{ \lambda} = 1$ the operator $\wcomp$ is explicitly given by
	\begin{equation}\label{defn:wcompLambda1}
	\left( \wcomp f \right) (z) = \psi_0\, e^{-\alpha \overline{a} \lambda z} f\big( a+ \lambda z\big),
	\end{equation}
	where we let $\psi_0 = \psi(0)$. 
	By comparing \eqref{defn:Tisometry} and \eqref{defn:wcompLambda1}, notice that 
	\begin{equation*}
	\wcomp = \psi_0 e^{\frac{\alpha}{2}\abs{a}^2} T_{\lambda, a}.
	\end{equation*}
	Hence the weighted composition operator $\wcomp$ is a constant multiple of an isometry and it follows that $\wcomp$ cannot be supercyclic.
\end{proof}

When $\abs{\lambda} = 1$, we recall for the Hilbert space $\Fh$ that Le~\cite[Theorem 3.3]{Le14} proved the weighted composition operator $\wcomp$ is a constant multiple of a unitary operator and hence normal.  
The following weakening of the definition of normality allows us to move beyond the Hilbertian case and consider general Banach spaces.
For a Banach space $X$, the bounded linear operator $U \colon X \to X $ is said to be \emph{paranormal} if
\begin{equation*}
\norm{ Ux }^2  \leq  \norm{ U^2x} \cdot \norm{x}
\end{equation*}
for all  $x \in X$. 
It is well known by results of Bourdon~\cite{Bou97} (cf.~\cite[p.~159]{GEP11}) that paranormal operators cannot be supercyclic.

An immediate consequence of Theorem \ref{thm:ConstIsometry} is that $\wcomp$ is a paranormal operator.

\begin{cor} \label{thm:wCompParanormal}
	Let $\wcomp$ be as in Theorem \ref{thm:ConstIsometry}. Then $\wcomp$ is a paranormal operator with the identity
	\begin{equation*}
	\norm{\wcomp f}^2 = \norm{\wcomp^2 f} \cdot \norm{f}
	\end{equation*}
	for all $f \in \F$, $1 \leq p \leq \infty$ and all $f \in \Fzero$.
\end{cor}

In general, in the Hilbert space setting, not every paranormal operator is normal. However, for $\wcomp$ acting on $\Fh$ with $\abs{\lambda} = 1$, it is a straightforward consequence of \cite[Theorem 3.3]{Le14} that $\wcomp$ is normal if and only if it is paranormal.

\section{Concluding remarks}

The property of \emph{weak supercyclicity} has been studied, for instance, in \cite{San04}, \cite{Shk07}, \cite[Chapter 10]{BM09} and \cite{KD18a}.  In particular, weak supercyclicity for (weighted) composition operators acting on various function spaces has been investigated in \cite{KHK10}, \cite{Bes13}, \cite{Bes14} and \cite{BJM20}. So a natural question arising from our investigation asks whether (weighted) composition operators can ever be weakly supercyclic on  Fock spaces? It would also be interesting to consider these questions in the higher dimensional setting of Fock spaces on $\C^n$ for $n \in \N$, as investigated by Carswell et al.~\cite{CMS03}.

In Section \ref{sec:paranormal} we show for $\abs{\lambda} = 1$ that weighted composition operators acting on Fock spaces are paranormal.  It is thus natural to ask, what is the complete characterisation of the paranormal weighted composition operators $\wcomp$ in the setting of Fock spaces?

\section*{Acknowledgements}

The authors with to thank Xavier Massaneda and Joaquim Ortega-Cerd\`{a} for insightful discussions and sharing their expert knowledge of Fock spaces.

%
%
%

\begin{thebibliography}{10}
	
	\bibitem{Ans95}
	S.~I. Ansari.
	\newblock Hypercyclic and cyclic vectors.
	\newblock {\em J. Funct. Anal.}, 128(2):374--383, 1995.
	
	\bibitem{AB97}
	S.~I. Ansari and P.~S. Bourdon.
	\newblock Some properties of cyclic operators.
	\newblock {\em Acta Sci. Math. (Szeged)}, 63(1-2):195--207, 1997.
	
	\bibitem{AL04}
	R.~Aron and M.~Lindstr\"{o}m.
	\newblock Spectra of weighted composition operators on weighted {B}anach spaces
	of analytic functions.
	\newblock {\em Israel J. Math.}, 141:263--276, 2004.
	
	\bibitem{BM09}
	F.~Bayart and {\'E}.~Matheron.
	\newblock {\em Dynamics of linear operators}, volume 179 of {\em Cambridge
		Tracts in Mathematics}.
	\newblock Cambridge University Press, Cambridge, 2009.
	
	\bibitem{Bel20}
	M.~J. Beltr\'{a}n-Meneu.
	\newblock Dynamics of weighted composition operators on weighted {B}anach
	spaces of entire functions.
	\newblock {\em J. Math. Anal. Appl.}, 492(1):124422, 2020.
	
	\bibitem{BGJJ16a}
	M.~J. Beltr\'{a}n-Meneu, M.~C. G\'{o}mez-Collado, E.~Jord\'{a}, and D.~Jornet.
	\newblock Mean ergodic composition operators on {B}anach spaces of holomorphic
	functions.
	\newblock {\em J. Funct. Anal.}, 270(12):4369--4385, 2016.
	
	\bibitem{BGJJ16b}
	M.~J. Beltr\'{a}n-Meneu, M.~C. G\'{o}mez-Collado, E.~Jord\'{a}, and D.~Jornet.
	\newblock Mean ergodicity of weighted composition operators on spaces of
	holomorphic functions.
	\newblock {\em J. Math. Anal. Appl.}, 444(2):1640--1651, 2016.
	
	\bibitem{BJM20}
	M.~J. Beltr\'{a}n-Meneu, E.~Jord\'{a}, and M.~Murillo-Arcila.
	\newblock Supercyclicity of weighted composition operators on spaces of
	continuous functions.
	\newblock {\em Collect. Math.}, 71(3):493--509, 2020.
	
	\bibitem{BCK10}
	C.~B\'{e}n\'{e}teau, B.~J. Carswell, and S.~Kouchekian.
	\newblock Extremal problems in the {F}ock space.
	\newblock {\em Comput. Methods Funct. Theory}, 10(1):189--206, 2010.
	
	\bibitem{BM95}
	L.~Bernal~Gonz\'{a}lez and A.~Montes-Rodr\'{\i}guez.
	\newblock Universal functions for composition operators.
	\newblock {\em Complex Variables Theory Appl.}, 27(1):47--56, 1995.
	
	\bibitem{BBMP15}
	N.~C. Bernardes, Jr., A.~Bonilla, V.~M{\"u}ller, and A.~Peris.
	\newblock Li-{Y}orke chaos in linear dynamics.
	\newblock {\em Ergodic Theory Dynam. Systems}, 35(6):1723--1745, 2015.
	
	\bibitem{Bes13}
	J.~B\`{e}s.
	\newblock Dynamics of composition operators with holomorphic symbol.
	\newblock {\em Rev. R. Acad. Cienc. Exactas F\'{\i}s. Nat. Ser. A Mat. RACSAM},
	107(2):437--449, 2013.
	
	\bibitem{Bes14}
	J.~B{\`e}s.
	\newblock Dynamics of weighted composition operators.
	\newblock {\em Complex Anal. Oper. Theory}, 8(1):159--176, 2014.
	
	\bibitem{BD11}
	J.~Bonet and P.~Doma\'{n}ski.
	\newblock A note on mean ergodic composition operators on spaces of holomorphic
	functions.
	\newblock {\em Rev. R. Acad. Cienc. Exactas F\'{\i}s. Nat. Ser. A Mat. RACSAM},
	105(2):389--396, 2011.
	
	\bibitem{BLW12}
	J.~Bonet, M.~Lindstr\"{o}m, and E.~Wolf.
	\newblock Norm-attaining weighted composition operators on weighted {B}anach
	spaces of analytic functions.
	\newblock {\em Arch. Math. (Basel)}, 99(6):537--546, 2012.
	
	\bibitem{Bou97}
	P.~S. Bourdon.
	\newblock Orbits of hyponormal operators.
	\newblock {\em Michigan Math. J.}, 44(2):345--353, 1997.
	
	\bibitem{BS90}
	P.~S. Bourdon and J.~H. Shapiro.
	\newblock Cyclic composition operators on {$H^2$}.
	\newblock In {\em Operator theory: operator algebras and applications, {P}art 2
		({D}urham, {NH}, 1988)}, volume~51 of {\em Proc. Sympos. Pure Math.}, pages
	43--53. Amer. Math. Soc., Providence, RI, 1990.
	
	\bibitem{BS97}
	P.~S. Bourdon and J.~H. Shapiro.
	\newblock Cyclic phenomena for composition operators.
	\newblock {\em Mem. Amer. Math. Soc.}, 125(596):x+105, 1997.
	
	\bibitem{CMS03}
	B.~J. Carswell, B.~D. MacCluer, and A.~Schuster.
	\newblock Composition operators on the {F}ock space.
	\newblock {\em Acta Sci. Math. (Szeged)}, 69(3-4):871--887, 2003.
	
	\bibitem{CH01}
	M.~D. Contreras and A.~G. Hern\'{a}ndez-D\'{\i}az.
	\newblock Weighted composition operators on {H}ardy spaces.
	\newblock {\em J. Math. Anal. Appl.}, 263(1):224--233, 2001.
	
	\bibitem{Cow83}
	C.~C. Cowen.
	\newblock Composition operators on {$H^{2}$}.
	\newblock {\em J. Operator Theory}, 9(1):77--106, 1983.
	
	\bibitem{CM95}
	C.~C. Cowen and B.~D. MacCluer.
	\newblock {\em Composition operators on spaces of analytic functions}.
	\newblock Studies in Advanced Mathematics. CRC Press, Boca Raton, FL, 1995.
	
	\bibitem{GGMR04}
	E.~A. Gallardo-Guti\'{e}rrez and A.~Montes-Rodr\'{\i}guez.
	\newblock The role of the spectrum in the cyclic behavior of composition
	operators.
	\newblock {\em Mem. Amer. Math. Soc.}, 167(791):x+81, 2004.
	
	\bibitem{GS16}
	E.~A. Gallardo-Guti\'{e}rrez and R.~Schroderus.
	\newblock The spectra of linear fractional composition operators on weighted
	{D}irichlet spaces.
	\newblock {\em J. Funct. Anal.}, 271(3):720--745, 2016.
	
	\bibitem{Gil20}
	C.~Gilmore.
	\newblock Linear dynamical systems.
	\newblock {\em Irish Math. Soc. Bull.}, 86:47--77, 2020.
	
	\bibitem{GEP11}
	K.-G. Grosse-Erdmann and A.~Peris~Manguillot.
	\newblock {\em Linear chaos}.
	\newblock Universitext. Springer, London, 2011.
	
	\bibitem{Gun07}
	G.~Gunatillake.
	\newblock Spectrum of a compact weighted composition operator.
	\newblock {\em Proc. Amer. Math. Soc.}, 135(2):461--467, 2007.
	
	\bibitem{Gun11}
	G.~Gunatillake.
	\newblock Invertible weighted composition operators.
	\newblock {\em J. Funct. Anal.}, 261(3):831--860, 2011.
	
	\bibitem{GI08}
	K.~Guo and K.~Izuchi.
	\newblock Composition operators on {F}ock type spaces.
	\newblock {\em Acta Sci. Math. (Szeged)}, 74(3-4):807--828, 2008.
	
	\bibitem{HK16}
	P.~V. Hai and L.~H. Khoi.
	\newblock Boundedness and compactness of weighted composition operators on
	{F}ock spaces {$\mathcal{F}^p(\mathbb{C})$}.
	\newblock {\em Acta Math. Vietnam.}, 41(3):531--537, 2016.
	
	\bibitem{HW73}
	H.~M. Hilden and L.~J. Wallen.
	\newblock Some cyclic and non-cyclic vectors of certain operators.
	\newblock {\em Indiana Univ. Math. J.}, 23:557--565, 1973/74.
	
	\bibitem{HLNS13}
	O.~Hyv\"{a}rinen, M.~Lindstr\"{o}m, I.~Nieminen, and E.~Saukko.
	\newblock Spectra of weighted composition operators with automorphic symbols.
	\newblock {\em J. Funct. Anal.}, 265(8):1749--1777, 2013.
	
	\bibitem{JPR87}
	S.~Janson, J.~Peetre, and R.~Rochberg.
	\newblock Hankel forms and the {F}ock space.
	\newblock {\em Rev. Mat. Iberoam.}, 3(1):61--138, 1987.
	
	\bibitem{JPZ17}
	L.~Jiang, G.~T. Prajitura, and R.~Zhao.
	\newblock Some characterizations for composition operators on the {F}ock
	spaces.
	\newblock {\em J. Math. Anal. Appl.}, 455(2):1204--1220, 2017.
	
	\bibitem{KHK10}
	Z.~Kamali, K.~Hedayatian, and B.~Khani~Robati.
	\newblock Non-weakly supercyclic weighted composition operators.
	\newblock {\em Abstr. Appl. Anal.}, page~14, 2010.
	\newblock Art. ID 143808.
	
	\bibitem{Kam75}
	H.~Kamowitz.
	\newblock The spectra of composition operators on {$H^{p}$}.
	\newblock {\em J. Funct. Anal.}, 18:132--150, 1975.
	
	\bibitem{Kam78}
	H.~Kamowitz.
	\newblock The spectra of a class of operators on the disc algebra.
	\newblock {\em Indiana Univ. Math. J.}, 27(4):581--610, 1978.
	
	\bibitem{KD18a}
	C.~S. Kubrusly and B.~P. Duggal.
	\newblock On weak supercyclicity {I}.
	\newblock {\em Math. Proc. R. Ir. Acad.}, 118A(2):47--63, 2018.
	
	\bibitem{Le14}
	T.~Le.
	\newblock Normal and isometric weighted composition operators on the {F}ock
	space.
	\newblock {\em Bull. Lond. Math. Soc.}, 46(4):847--856, 2014.
	
	\bibitem{MS20}
	T.~Mengestie and W.~Seyoum.
	\newblock Topological and dynamical properties of composition operators.
	\newblock {\em Complex Anal. Oper. Theory}, 14(1):paper no. 2, 27, 2020.
	
	\bibitem{NRW87}
	E.~Nordgren, P.~Rosenthal, and F.~S. Wintrobe.
	\newblock Invertible composition operators on {$H^p$}.
	\newblock {\em J. Funct. Anal.}, 73(2):324--344, 1987.
	
	\bibitem{Nor68}
	E.~A. Nordgren.
	\newblock Composition operators.
	\newblock {\em Canadian J. Math.}, 20:442--449, 1968.
	
	\bibitem{Rez11}
	H.~Rezaei.
	\newblock Chaotic property of weighted composition operators.
	\newblock {\em Bull. Korean Math. Soc.}, 48(6):1119--1124, 2011.
	
	\bibitem{San04}
	R.~Sanders.
	\newblock Weakly supercyclic operators.
	\newblock {\em J. Math. Anal. Appl.}, 292(1):148--159, 2004.
	
	\bibitem{SMB20}
	W.~Seyoum, T.~Mengestie, and J.~Bonet.
	\newblock Mean ergodic composition operators on generalized {F}ock spaces.
	\newblock {\em Rev. R. Acad. Cienc. Exactas F\'{\i}s. Nat. Ser. A Mat. RACSAM},
	114(1):paper no. 6, 11, 2020.
	
	\bibitem{Sha87}
	J.~H. Shapiro.
	\newblock The essential norm of a composition operator.
	\newblock {\em Ann. of Math. (2)}, 125(2):375--404, 1987.
	
	\bibitem{Sha93}
	J.~H. Shapiro.
	\newblock {\em Composition operators and classical function theory}.
	\newblock Universitext. Springer-Verlag, New York, 1993.
	
	\bibitem{Shk07}
	S.~Shkarin.
	\newblock Non-sequential weak supercyclicity and hypercyclicity.
	\newblock {\em J. Funct. Anal.}, 242(1):37--77, 2007.
	
	\bibitem{Tie20a}
	P.~T. Tien.
	\newblock Composition operators on weighted banach spaces of entire functions.
	\newblock {\em Complex Variables and Elliptic Equations}, in press:1--17, 2020.
	
	\bibitem{Tie20b}
	P.~T. Tien.
	\newblock The iterates of composition operators on {B}anach spaces of
	holomorphic functions.
	\newblock {\em J. Math. Anal. Appl.}, 487(1):123945, 23, 2020.
	
	\bibitem{TK19}
	P.~T. Tien and L.~H. Khoi.
	\newblock Weighted composition operators between different {F}ock spaces.
	\newblock {\em Potential Anal.}, 50(2):171--195, 2019.
	
	\bibitem{Uek07}
	S.-I. Ueki.
	\newblock Weighted composition operator on the {F}ock space.
	\newblock {\em Proc. Amer. Math. Soc.}, 135(5):1405--1410, 2007.
	
	\bibitem{Uek10}
	S.-I. Ueki.
	\newblock Weighted composition operators on some function spaces of entire
	functions.
	\newblock {\em Bull. Belg. Math. Soc. Simon Stevin}, 17(2):343--353, 2010.
	
	\bibitem{YR07}
	B.~Yousefi and H.~Rezaei.
	\newblock Hypercyclic property of weighted composition operators.
	\newblock {\em Proc. Amer. Math. Soc.}, 135(10):3263--3271, 2007.
	
	\bibitem{Zhu12}
	K.~Zhu.
	\newblock {\em Analysis on {F}ock spaces}, volume 263 of {\em Graduate Texts in
		Mathematics}.
	\newblock Springer, New York, 2012.
	
\end{thebibliography}

\end{document}